\def\VERSIONNUMBER{1.2}
\def\VERSION{Version \VERSIONNUMBER}
\def\DATE{28 September 1998, revised on 30 January 1999}
\def\EPRINT{{\tt math.QA/9809157}}
\documentclass[]{amsart}
\usepackage{amsmath,amssymb
}

\newcommand\commentout[1]{}
\parskip \medskipamount

\newcommand\ad{\operatorname{ad}}

\newcommand\borel{{\mathfrak b}}

\newcommand\bos{{{\mathrm b}{\mathrm o}{\mathrm s}}}
\newcommand\Bos{\operatorname{Bos}}
\newcommand\bosFock{{\mathcal F}^{\bos}}
\newcommand\bosket{\rangle^{\bos}}

\newcommand\bcdot{{\raise0.2ex \hbox{\bf.}}}
\newcommand\bdot{{\raise0.4ex \hbox{\bf.}}}
\newcommand\C{{\mathbb C}}
\newcommand\Complex{{\mathbb C}}
\newcommand\caff{{{\mathrm c}_{\mathrm a}}}
\newcommand\cvir{{c_{\mathrm V}}}
\newcommand\CB{\operatorname{CB}\nolimits}

\newcommand\ch{\operatorname{ch}\nolimits}
\newcommand\CC{\operatorname{CC}\nolimits}
\newcommand\Ctimes{{\mathbb C}^\times}
\newcommand\D{{\mathcal D}}

\newcommand\End{\operatorname{End}\nolimits}
\newcommand\eps{\varepsilon}
\newcommand\F{{\mathcal F}}

\newcommand\Fock{{\mathcal F}}

\newcommand\g{{\mathfrak g}}

\newcommand\germanS{{\mathfrak S}}
\newcommand\gh{{{\mathrm g}{\mathrm h}}}
\newcommand\Gh{\operatorname{Gh}\nolimits}

\newcommand\gHDpr{\g^{H,D}_{\prin}}    
\newcommand\ghFock{{\mathcal F}^{\gh}}
\newcommand\ghvac{|0\rangle^{\gh}}

\newcommand\h{{\mathfrak h}}
\newcommand\Heis{{\mathcal H}}

\newcommand\Hom{\operatorname{Hom}\nolimits}

\newcommand\id{{{\mathrm i}{\mathrm d}}}

\newcommand\Integer{{\mathbb Z}}

\newcommand\khat{{\hat k}} 
\newcommand\K{{\mathcal K}}

\newcommand\Line{{\mathcal L}}
\newcommand\m{{\mathfrak m}}

\ifx\mathit\undefined
\newcommand\mathit{} 
\fi
\newcommand\n{{\mathfrak n}}
\newcommand\Natural{{\mathbb N}}

\newcommand\np{{:}} 
\newcommand\NP{{\genfrac{}{}{0pt}{}{\circ}{\circ}}} 
\renewcommand\O{{\mathcal O}}

\newcommand\onto{\twoheadrightarrow}

\newcommand\calP{{\mathcal P}}
\ifx\pd\undefined
\newcommand\pd[2]{\frac{\partial #1}{\partial #2}}
\else
\renewcommand\pd[2]{\frac{\partial #1}{\partial #2}} 
\fi

\newcommand\pr{\pi}
\newcommand\prin{{\dot X}} 
\newcommand\Res{\mathop{\operatorname{Res}}\nolimits}
\newcommand\Scr{Scr}
\newcommand\scr{scr}
\renewcommand\setminus{\smallsetminus}
\newcommand\simeqq{\cong}
\newcommand\Span{\mathop{\operatorname{Span}}\nolimits}
\newcommand\Suga{{\mathrm S}ug}

\newcommand\tensor{\otimes}

\newcommand\Tr{\operatorname{Tr}\nolimits}

\newcommand\Vir{{\mathrm V}{\mathrm i}{\mathrm r}}
\newcommand\Waki{{\mathrm W}ak}

\newcommand\X{{\mathfrak X}}
\newcommand\Xtilde{{\widetilde{\mathfrak X}}}
\newcommand\Z{{\mathbb Z}}

\newtheorem{thm}{Theorem}[section]
\newtheorem{lem}[thm]{Lemma}
\newtheorem{prop}[thm]{Proposition}
\newtheorem{cor}[thm]{Corollary}

\theoremstyle{definition}
\newtheorem{defn}[thm]{Definition}

\theoremstyle{remark}
\newtheorem{example}[thm]{Example}
\newtheorem{rem}[thm]{Remark}

\numberwithin{equation}{section}

\newcommand\thmref[1]{Theorem~\ref{#1}}
\newcommand\propref[1]{Proposition~\ref{#1}}
\newcommand\secref[1]{Section~\ref{#1}}
\newcommand\corref[1]{Corollary~\ref{#1}}
\newcommand\lemref[1]{Lemma~\ref{#1}}
\newcommand\defref[1]{Definition~\ref{#1}}

\newcommand\appref[1]{Appendix~\ref{#1}}

\begin{document}

\title[Bosonization and KZB equations]
{Bosonization and integral representation\\
of solutions of the Knizhnik-Zamolodchikov-Bernard equations}


\author{Gen KUROKI} %
\address{Mathematical Institute, Tohoku University, Sendai 980, JAPAN}
\author{Takashi TAKEBE}%
\address{Department of Mathematical Sciences,
the University of Tokyo, Komaba, Tokyo 153-8914, JAPAN}

\date{\DATE{}\quad
(\VERSION, \EPRINT, University of Tokyo preprint UTMS 98-39)}
\ifx\ClassWarning\undefined
\maketitle 
\fi

\begin{abstract}
We construct an integral representation of solutions of the
Knizhnik-Zamolodchikov-Bernard equations, using the Wakimoto modules.
\end{abstract}

\ifx\ClassWarning\undefined\else
\maketitle
\fi

\centerline{\VERSION.}


\tableofcontents

\parskip = 0pt
\begin{small}
\tableofcontents
\end{small}
\parskip = \medskipamount
\bigskip

\setcounter{section}{-1}
\section{Introduction}
\label{intro}
The purpose of this paper is to construct an integral representation of
solutions of the Knizhnik-Zamolodchikov-Bernard equations, using the
Wakimoto modules.

Correlation functions of the (chiral) Wess-Zumino-Witten models satisfy a
system of differential equations. In the genus zero case, it is the
well-known Knizhnik-Zamolodchikov (KZ) equations
\cite{kni-zam:84}, \cite{tsu-kan:88}. Bernard \cite{ber:88-1} found a
system of equations for
the genus one case, which is now called the Knizhnik-Zamolodchikov-Bernard
(KZB) equations. In general, it is known (cf.\ \cite{tu-ue-ya:89}) that
correlation functions are solutions of a holonomic system over the base
space of a family of Riemann surfaces with marked points and principal
bundles on them. See also \cite{ber:88-2} and \cite{fel:98}.

There are a vast amount of works on the KZ equations, among which are
studies on integral representation of solutions. There are several
different approaches to this subject. One is from the viewpoint of the
theory of hypergeometric type integrals; e.g., \cite{djmm:89},
\cite{mat:90}, \cite{sch-var:90}, \cite{sch-var:91}. Another approach
comes from free field theories on Riemann surfaces and the Wakimoto
realization of affine Lie algebras; e.g., \cite{mar:89}, \cite{gmmos:90},
\cite{ber-fel:90},
\cite{a-t-y:91}, \cite{awa:91}. The third one is the off-shell Bethe
Ansatz developed in \cite{bab-flu:94}, the representation theoretical
meaning of which was clarified by \cite{f-f-r:94}. See also
\cite{resh-var:95} and \cite{cher:91}. 

The first approach in the genus one case was pursued by Felder and
Varchenko in \cite{fel-var:95}, while Babujian et al.\ \cite{b-ls-p:98}
study the off-shell Bethe Ansatz approach of the KZB equations (with an
additional term). Our goal in this paper is to apply the second approach
to the genus one case to obtain an integral representation of solutions of
the KZB equations. In the genus zero case, an integral over a suitable
(twisted) cycle of a matrix element of product of vertex operators and
screening charges gives a solution of the KZ equations. In the genus one
case, we take a twisted trace instead of a matrix element to give a
solution of the KZB equation in the integral form. We apply the method of
the screening current Ward identity used in \cite{a-t-y:91} and
\cite{awa:91} mutatis mutandis, and obtain an explicit formula for this
integral representation. For the $sl(2)$ case Bernard and Felder found the 
same result in \cite{ber-fel:90} by using Wakimoto modules.

The paper is organized as follows.  In \secref{wakimoto}, mainly following
\cite{kur:91}, we review several fundamental techniques in the conformal
field theory, especially the free field realization of the affine Lie
algebra found by Wakimoto \cite{wak:86}, Feigin and Frenkel
\cite{fei-fre:88}, \cite{fei-fre:90}. We state the problem in
\secref{wzw}. Namely we formulate the Wess-Zumino-Witten model on elliptic
curves, following \cite{fel-wie:96}, and give a definition of $N$-point
functions. The KZB equations are introduced as a system of equations
satisfied by $N$-point functions. \secref{N-point-function:wakimoto} is
the main part of this paper, where we give an integral representation of
an $N$-point function on elliptic curves
(\thmref{int-rep-of-N-point-function}). If we restrict an $N$-point
function to a certain submodule, it is a solution of the KZB equations. We
thus give an integral representation of solutions of the KZB equations
(\thmref{integral-rep}). To write down the integrand explicitly, we use
the screening current Ward identity on elliptic curves. \appref{theta} is
a table of theta functions used in the paper. In \appref{coherent-states}
we review the method of coherent states well known in the string
theory. We also compute the one-loop correlation function of vertex
operators.

\section{Bosonization and Wakimoto modules}
\label{wakimoto}

In this section we review basic facts about the Wakimoto representations of
the affine Lie algebras, following \cite{kur:91}. See also
\cite{wak:86}, \cite{fei-fre:88}, \cite{fei-fre:90}, \cite{f-f-r:94}.

\subsection{Notations for the finite dimensional algebra}
\label{notations:fin-dim}

Here we recall fundamental facts about finite dimensional simple Lie
algebras to fix the notations.

Let $\g$ be a finite dimensional simple Lie algebra of rank $l$, $\h$
its Cartan subalgebra and
\begin{equation}
    \g = \h \oplus \bigoplus_{\alpha \in \Delta} \g_\alpha
\label{root-decomp}
\end{equation}
the root space decomposition, where $\Delta$ is the set of roots. The
Cartan-Killing form is denoted by $(\cdot | \cdot )$, through which we
identify $\h$ and its dual space $\h^\ast$.  We fix the simple roots
$\{\alpha_1,\ldots,\alpha_l\}$, Chevalley generators $\{H_i, E_i,
F_i\}_{i=1,\dots,l}$ and a basis $e_\alpha$ of $\g_\alpha$, such that
$e_{\alpha_i} = E_i$ for $i=1,\ldots,l$ and
$(e_\alpha|e_{-\alpha'})=\delta_{\alpha,\alpha'}$. The set of positive and
negative roots are denoted by $\Delta_+ = \{\beta_1, \dots, \beta_s \}$
and $\Delta_-$ respectively. The Borel and nilpotent subalgebras
corresponding to $\Delta_\pm$ are denoted by $\borel_\pm$, $\n_\pm$ as
usual.

Let $G$, $B_\pm$ and $N_\pm$ be an algebraic group corresponding to $\g$,
the subgroups corresponding to $\borel_\pm$ and to $\n_\pm$. As is well
known, there exists a Lie algebra homomorphism $R_\lambda$ from $\g$ to
the sheaf of twisted differential operators $\D_\lambda$ on the flag
variety $B_-\backslash G$ once one fixes a dual vector
$\lambda\in\h^\ast$. Denote the base point $[B_+]\in B_-\backslash G$ by
$o$. With respect to the coordinate on $o N_+ \cong \n_+$, a big cell of
$B_-\backslash G$, introduced by the exponential map:
$$
    \Complex^{\Delta_+} \owns
    (x^\alpha)_{\alpha\in\Delta_+} \mapsto 
    o \exp(x^{\beta_1} e_{\beta_1}) \cdots \exp(x^{\beta_s} e_{\beta_s})
    \in o N_+,
$$
the twisted differential operator $R_\lambda(X)$ for $X \in \g$ is
represented by a first order differential operator acting on the space of
polynomials $\Complex[x^\alpha; \alpha \in \Delta_+]$:
\begin{equation}
    R_\lambda(X) = R(X; x, \partial_x, \lambda),
\label{diff-op-rep}
\end{equation}
where $\partial_x = (\partial/\partial x^\alpha)_{\alpha\in\Delta_+}$. The
operator $R(X;x,\partial_x,\lambda)$ is a polynomial in $X$, $x$,
$\partial_x$ and $\lambda(H_i)$ ($i=1,\ldots,l$). More explicitly, there
are polynomials $R_\alpha(X;x)$ in $x$ for $X\in\g$, $\alpha\in\Delta_+$
such that
\begin{equation}
\begin{aligned}
    R_\lambda(E_i) 
    &= \sum_{\alpha\in\Delta_+}
       R_\alpha(E_i;x) \frac{\partial}{\partial x_\alpha},
\\
    R_\lambda(F_i) 
    &= \sum_{\alpha\in\Delta_+}
       R_\alpha(F_i;x) \frac{\partial}{\partial x_\alpha}
    + x_{\alpha_i} \lambda(H_i),
\\
    R_\lambda(H) 
    &= - \sum_{\alpha\in\Delta_+}
         \alpha(H) x_\alpha \frac{\partial}{\partial x_\alpha}
    + \lambda(H),
\end{aligned}
\label{explicit-diff-op}
\end{equation}
for Chevalley generators $E_i$, $F_i$ ($i= 1,\ldots,l$) and $H\in\h$.
Note that $R_\lambda(E_i)$ does not depend on $\lambda$. Hence we
sometimes omit the suffix and denote it by $R(E_i)$.

The nilpotent subgroup $N_+$ acts on the big cell from the left as
$$
    n \cdot (o a) = ona \qquad \text{for }n,a\in N_+.
$$
The infinitesimal action of $\n_+$ induced by this action is denoted by
$\Scr$:
\begin{equation}
    n_+ \owns X \mapsto 
    \Scr(X;x,\partial_x) \in \Complex[x,\partial_x].
\label{nilp-left-action}
\end{equation}

\subsection{Ghosts and free bosons}
\label{ghost-boson}

Let us introduce the algebra of (bosonic) ghosts, $\widehat\Gh(\g)$, and
the algebra of free bosons, $\widehat\Bos(\g)$.

The generators of the algebra $\widehat\Gh(\g)$ are $\beta_\alpha[m]$ and
$\gamma^\alpha[m]$ ($\alpha\in\Delta_+$, $m\in\Integer$) satisfying the
canonical commutation relations:
\begin{equation}
    [\beta_\alpha[m], \gamma^{\alpha'}[n]] 
    = \delta_\alpha^{\alpha'} \delta_{m+n,0} \cdot 1,
\label{ghost-relation}
\end{equation}
for $\alpha, \alpha' \in \Delta_+$ and $m,n\in\Integer$. We call the formal 
generating functions of generators,
\begin{equation}
    \beta_\alpha(z) = \sum_{m\in\Integer} z^{-m-1}\beta_\alpha[m], \qquad
    \gamma^\alpha(z) = \sum_{m\in\Integer} z^{-m}\gamma^\alpha[m],
\label{def:beta-gamma}
\end{equation}
{\em ghost fields}. They satisfy the following operator product
expansions:
\begin{equation}
    \beta_\alpha(z) \gamma^{\alpha'}(w) 
    \sim \frac{\delta_\alpha^{\alpha'}}{z-w}.
\label{ope:ghost}
\end{equation}
The {\em ghost Fock space} $\ghFock$ is defined as a left
$\widehat\Gh(\g)$-module generated by the vacuum vector $\ghvac$,
satisfying
\begin{equation}
    \beta_\alpha[m]  \ghvac = 0, \qquad 
    \gamma^\alpha[n] \ghvac = 0
\label{def:ghost-vac}
\end{equation}
for any $\alpha\in\Delta_+$, $m\geqq 0$, $n>0$.

The algebra $\widehat\Bos(\g)$ is generated by $\phi_i[m]$
($i=1,\dots,l$, $m\in\Integer$), the defining relation of which is
\begin{equation}
    [\phi_i[m], \phi_j[n]] = \kappa (H_i|H_j) m \delta_{m+n,0}\cdot 1,
\label{boson-relation}
\end{equation}
where $\kappa$ is a non-zero complex parameter. We extend the algebra to
$\widetilde\Bos(\g)$ by adding the boost operator $e^{p_i}$ and its
logarithm $p_i$ which satisfies the relation:
\begin{equation}
    [\phi_i[m], e^{p_j}] = \kappa (H_i|H_j) \delta_{m,  0} e^{p_j},
    \qquad 
    [\phi_i[m],    p_j ] = \kappa (H_i|H_j) \delta_{m,  0}.
\label{boost-relation}
\end{equation}

Fields
\begin{equation}
\begin{aligned}
    \phi_i(z) &:= \kappa p_i + \phi_i[0] \log z +
    \sum_{m\in\Integer\setminus\{0\}} \frac{z^{-m}}{-m} \phi_i[m],\\
    \partial\phi_i(z) &:= \sum_{m\in\Integer} z^{-m-1} \phi_i[m]
\end{aligned}
\label{def:boson}
\end{equation}
are important generating functions of generators of this algebra. The
field $\phi_i(z)$ is called the {\em free boson field\/}. For any $H=
\sum_{i=1}^l a_i H_i \in\h$, we also use the notations like
$$
    \phi[H;m] = \sum_{i=1}^l a_i \phi_i[m], \qquad
    p[H]      = \sum_{i=1}^l a_i p_i,       \qquad
    \phi(H;z) = \sum_{i=1}^l a_i \phi_i(z).
$$
The free boson fields satisfy the following operator product expansion.
\begin{equation}
\begin{aligned}
    \phi(H;z) \phi(H';w) &\sim \kappa (H|H') \log(z-w),\\
    \partial\phi(H;z) \partial\phi(H';w) 
    &\sim \frac{\kappa (H|H')}{(z-w)^2},
\end{aligned}
\label{ope:boson}
\end{equation}
for any $H, H' \in \h$. For a dual vector $\lambda\in\h^\ast$, the {\em
boson Fock space} $\bosFock_\lambda$ with momentum $\lambda$ is defined as
a left $\widehat\Bos(\g)$-module generated by the vacuum vector $|\lambda\bosket$,
satisfying
\begin{equation}
    \phi_i[m] |\lambda\bosket = 0, \qquad
    \phi_i[0] |\lambda\bosket = \lambda(H_i) |\lambda\bosket
\label{def:boson-vac}
\end{equation}
for any $i=1,\dots,l$, $m>0$. The boost operator $e^{p_i}$ acts on the
direct sum of Fock spaces $\bigoplus_{\lambda} \bosFock_\lambda$ by
shifting the momentum:
\begin{equation}
    e^{p_i} |\lambda\bosket = |\lambda+H_i\bosket.
\label{def:boost-vac}
\end{equation}

The {\em normal ordered product} $\np P \np$ of a monomial $P$ of
$\beta_\alpha[m]$'s, $\gamma^\alpha[m]$'s, $\phi_i[m]$'s and $e^{p_i}$'s
is defined by putting annihilation operators of $\ghvac$ and
$|\lambda\bosket$ ($\beta_\alpha[m]$ ($m\geqq0$), $\gamma^\alpha[m]$
($m>0$), $\phi_i[m]$ ($m>0$)) and $\phi_i[0]$ appearing in $P$ to the
right side in the product.

For example, the {\em bosonic vertex operator} is defined by
\begin{equation}
\begin{split}
    V(\lambda;z) &:= \np e^{\frac{\phi(\lambda;z)}{\kappa}} \np \\
    &=
    \exp \left(
         \frac{1}{\kappa}\sum_{m<0}\frac{z^{-m}}{-m} \phi[\lambda;m]
         \right)
    e^{p[\lambda]} z^{\frac{1}{\kappa} \phi[\lambda;0]}
    \exp \left(
         \frac{1}{\kappa}\sum_{m>0}\frac{z^{-m}}{-m} \phi[\lambda;m]
         \right).
\end{split}
\label{def:boson-vo}
\end{equation}
We introduce the following notations for later use.
\begin{align}
    V(\lambda;z) &= \tilde V(\lambda;z) V_0(\lambda;z),
\label{boson-vo:decomposition}
\\
    \tilde V(\lambda;z) 
    &:= \np e^{\frac{\tilde\phi(\lambda;z)}{\kappa}} \np
    =
    \exp \left(
         \frac{1}{\kappa}\sum_{m<0}\frac{z^{-m}}{-m} \phi[\lambda;m]
         \right)
    \exp \left(
         \frac{1}{\kappa}\sum_{m>0}\frac{z^{-m}}{-m} \phi[\lambda;m]
         \right),
\label{def:tildeV}
\\
    V_0(\lambda;z) &:=
    e^{p[\lambda]} z^{\frac{1}{\kappa} \phi[\lambda;0]},
\label{def:V0}
\end{align}
where $\tilde\phi(\lambda;z) = \phi(\lambda;z) - \kappa p[\lambda] -
\phi[\lambda;0]$ is the non-zero mode part of $\phi(\lambda;z)$. 

\subsection{Bosonization and Wakimoto modules}
\label{bosonization}

Bosonizing the differential operators $R(X;x,\partial_x,\lambda)$ in
\secref{notations:fin-dim} by ghosts and free bosons in
\secref{ghost-boson} gives the Wakimoto realization of the affine Lie
algebra $\hat\g$. 

Define {\em current operator} $X(z)$ ($X\in\g$) and the {\em
energy-momentum tensor} $T(z)$ by
\begin{align}
    X(z) &:= \np R(X;\gamma(z), \beta(z), \partial\phi(z)) \np \qquad
    \text{for }X = E_i, H_i \text{ and }i=1,\dots,l,
\label{def:X(z)}
\\
    F_i(z)&:= \np R(F_i;\gamma(z), \beta(z), \partial\phi(z)) \np
            + c_i \partial\gamma^{\alpha_i}(z)
    \qquad \text{for }i=1,\dots,l,
\label{def:F-i(z)}
\\
    T(z) &:= T^{\gh}(z) + T^\phi(z),
\label{def:T(z)}
\\
    T^\gh(z)&:=
    \sum_{\alpha\in\Delta_+} 
    \np \partial \gamma^\alpha(z) \beta_\alpha(z) \np,
\label{def:Tgh(z)}
\\
    T^\phi(z) &:=
    \frac{1}{2\kappa}
    \sum_{i=1}^l \np \partial\phi(H_i;z) \partial\phi(H^i;z)\np 
    -
    \frac{1}{2\kappa}
    \partial^2\phi(2\rho;z),
\label{def:Tphi(z)}
\end{align}
where $\{H^i\}_{i=1,\ldots,l}$ is a basis of $\h$ dual to $\{H_i\}$ with
respect to $(\cdot | \cdot)$, $\rho$ is the half sum of positive roots
($\h$ and $\h^\ast$ are identified via the inner product) and
$\{c_i\}_{i=1,\dots,l}$ is a set of constants to be determined. More
explicitly, we have from \eqref{explicit-diff-op}
\begin{align}
    E_i(z) 
    &= \sum_{\alpha\in\Delta_+}
       \np R_\alpha(E_i;\gamma(z)) \beta_\alpha(z) \np,
\label{explicit:E-i(z)}
\\
    F_i(z) 
    &= \sum_{\alpha\in\Delta_+}
       \np R_\alpha(F_i;\gamma(z)) \beta_\alpha(z) \np
    + \gamma_{\alpha_i}(z) \partial\phi_i(z)
    + c_i \partial\gamma^{\alpha_i}(z),
\label{explicit:F-i(z)}
\\
    H(z)
    &= H^\gh(z)
    + \partial\phi(H;z), \qquad
    H^\gh(z):= - \sum_{\alpha\in\Delta_+}
               \alpha(H) \np \gamma^\alpha(z) \beta_\alpha(z) \np
\label{explicit:H(z)}
\end{align}
for Chevalley generators $E_i$, $F_i$ ($i= 1,\ldots,l$) and $H\in\h$.

We expand these series in the following way: for $X\in\g$,
\begin{equation}
    X(z) = \sum_{m\in\Integer} z^{-m-1} X[m], \qquad
    T(z) = \sum_{m\in\Integer} z^{-m-2} T[m].
\label{mode-expansion}
\end{equation}
The coefficients $X[m]$ and $T[m]$ belong to a certain completion of the
algebra $\widehat\Gh(\g) \tensor \widehat\Bos(\g)$.

\begin{thm}{%
\cite{wak:86}, \cite{fei-fre:88}, \cite{fei-fre:90}, \cite{kur:91}.}
There exists a unique set of constants $\{c_i\}_{i=1,\dots,l}$, such that
a Lie algebra homomorphism from the affine Lie algebra
$\hat\g=\g\tensor\Complex[t,t^{-1}]\oplus \Complex\hat k$ to a completion
of $\widehat\Gh(\g) \tensor \widehat\Bos(\g)$ can be defined by
\begin{equation}
    \omega(X\tensor t^m) = X[m], \qquad
    \omega(\hat k) = \kappa - h^\vee,
\label{wakimoto-hom-km}
\end{equation}
for all $X\in\g$, $m\in\Integer$, where $\hat k$ is the center of $\hat\g$
and $h^\vee$ is the dual Coxeter number of $\g$. Moreover, the
energy-momentum tensor $T(z)$ defined by \eqref{def:T(z)} coincides with
the image of $T_{\Suga}(z)$ in $U\hat\g$ defined by the Sugawara
construction:
\begin{equation}
    T(z) = \omega(T_{\Suga}(z)), \qquad
    T_{\Suga}(z) 
    := \frac{1}{2\kappa}\sum_{p=1}^{\dim \g} \NP J_p(z) J^p(z)\NP,
\label{em-tensor}
\end{equation}
and $T[m]$'s generate the Virasoro algebra $\Vir$ with the central charge
$\cvir=\dim\g - 12(\rho|\rho)/\kappa = k\dim\g/\kappa$. Here $J_p(z) =
\sum_{m\in\Integer} z^{-m-1}J_p\tensor t^m$, $\{J_p\}$ is a basis of $\g$,
$\{J^p\}$ is its dual basis with respect to $(\cdot|\cdot)$ and the symbol
$\NP \ \NP$ is the normal ordered product in $U\hat\g$.
\end{thm}

Namely, $\omega$ can be extended to a Lie algebra homomorphism from
$\hat\g\oplus\Vir$ to a completion of
$\widehat\Gh(\g)\oplus\widehat\Bos(\g)$ such that
\begin{equation}
    T[m] = \omega(T_{\Suga}[m]).
\label{wakimoto-hom-vir}
\end{equation}

Therefore Kac-Moody current operators satisfy the operator product
expansions:
\begin{equation}
    X(z) Y(w) \sim
    \frac{k (X|Y)}{(z-w)^2} + \frac{[X,Y](w)}{z-w},
\label{ope:KM-current}
\end{equation}
where $X, Y \in \g$ and $k=\kappa - h^\vee$, and the energy-momentum tensor
satisfies
\begin{align}
    T(z) T(w) &\sim 
      \frac{\cvir/2}{(z-w)^4} 
    + \frac{2T(w)}{(z-w)^2}
    + \frac{\partial T(w)}{z-w},
\label{ope:em}
\\
    T(z) X(w) &\sim
    \frac{X(w)}{(z-w)^2} + \frac{\partial X(w)}{z-w}.
\label{ope:em-KM}
\end{align}
We can regard $\ghFock \tensor \bosFock_\lambda$ as a representation of
$\hat\g$ of level $k=\kappa-h^\vee$ through $\omega$.
\begin{defn}
\label{def:wakimoto-rep}
Denote $\ghFock \tensor \bosFock_\lambda$ by $\Waki_{\lambda,k}$ and call
it a {\em Wakimoto module} of level $k$, weight $\lambda$.
\end{defn}

There is a $\g$-submodule generated by $\ghvac\tensor|\lambda\bosket$,
which is spanned by $\prod_{\alpha\in\Delta_+}
\gamma^\alpha[0]^{I(\alpha)}\ghvac\tensor|\lambda\bosket$
($I(\alpha)\in\Natural$), and is isomorphic to the dual Verma module
$M^\ast_\lambda$ of $\g$ with the highest weight $\lambda$. (See
Proposition 4.4 of \cite{kur:91}.) We denote it by
$\Waki^{0}_{\lambda,k}$:
\begin{equation}
    \Waki^{0}_{\lambda,k} := 
    \Span_\Complex\Bigl\{
    \prod_{\alpha\in\Delta_+}
    \gamma^\alpha[0]^{I(\alpha)}\ghvac\tensor|\lambda\bosket \,\Big|\,
    (I(\alpha))_{\alpha\in\Delta_+}\in \Natural^{\Delta_+}
    \Bigr\}
    \simeqq M^\ast_\lambda.
\label{def:wakimoto-finite}
\end{equation}
It is easy to show that for any $m>0$ and $X\in\g$,
\begin{equation}
    X[m] \Waki^{0}_{\lambda,k} = 0,
\label{nilp=zero:wakimoto}
\end{equation}
and the quadratic Casimir operator $C_2 = \sum_p J_p J^p$ acts as a
multiplication operator:
\begin{equation}
    C_2|_{\Waki^{0}_{\lambda,k}}
    =
    (\lambda|\lambda+2\rho) \id.
\label{casimir-on-wakimoto}
\end{equation}

\subsection{State-operator correspondence}
\label{state-operator-correspondence}

Let us recall the state-operator correspondence in the two dimensional
conformal field theories. A primary field generates a highest weight
representation of the algebra of symmetries of the theory in a space of
operator valued functions (``local operators''). In our case, the algebra
of symmetries is $\widehat\Gh(\g) \oplus \widehat\Bos(\g)$ and the
representation space is:
\begin{equation}
\begin{split}
    \O_\lambda&:=
    \Span_\Complex\{
    \widehat{x_1[m_1]}\cdots\widehat{x_n[m_n]} V(\lambda;z) |
    n\in\Natural;
    x_i=\beta_\alpha,\gamma^\alpha \text{ or }\phi_j 
    \text{ for certain }\alpha, j;
    m_i\in\Integer
    \} \\
    &\subset \widehat\Gh(\g) \oplus \widehat\Bos(\g)((z)),
\end{split}
\label{def:O-lambda}
\end{equation}
where $V(\lambda;z)$ is defined by \eqref{def:boson-vo} and the action of
$x[m]$, denoted by $\widehat{x[m]}$, is defined by
\begin{equation}
\begin{split}
    \widehat{x[m]} \Phi(z)
    &:=
    \Res_{\zeta=z}(\zeta-z)^{m+h-1} x(\zeta) \Phi(z), \qquad
    (x=\beta_\alpha, \gamma^\alpha),
\\
    \widehat{\phi_i[m]} \Phi(z)
    &:=
    \Res_{\zeta=z}(\zeta-z)^{m} \partial\phi_i(\zeta) \Phi(z),
\end{split}
\label{def:operator-realization}
\end{equation}
where $h$ is the conformal spin of the field $x(z) = \sum_{n\in \Integer}
x[n] z^{-n-h}$ (cf.\ \eqref{def:beta-gamma}) and $\Phi(z) \in \O_\lambda$.
An element of $\O_\lambda$ maps $\Waki_{\mu,k}$ to $\Waki_{\lambda+\mu,k}$
for any $\mu$:
\begin{equation}
    \O_\lambda 
    \subset \Hom_\Complex(\Waki_{\mu,k},\Waki_{\lambda+\mu,k}).
\label{O-lambda-act-on-wakimoto}
\end{equation}
Because of the operator product expansions, \eqref{ope:ghost},
\eqref{ope:boson}, $\O_\lambda$ is a $\widehat\Gh(\g) \oplus
\widehat\Bos(\g)$ module:
\begin{lem}
\label{operator-realization}
(i) For any fields $x, y = \beta_\alpha, \gamma^\alpha, \phi_i$ and any
integers $m,n\in\Integer$, we have
$$
    [\widehat{x[m]}, \widehat{y[n]}]
    =
    [x[m],y[n]]^{\wedge} \in \End_\Complex(\O_\lambda).
$$

(ii)
$\O_\lambda$ is generated by $V(\lambda;z)$, which satisfies
\begin{equation}
    \widehat{\beta_\alpha [m]} V(\lambda;z) =
    \widehat{\gamma^\alpha[n]} V(\lambda;z) =
    \widehat{\phi_i       [m]} V(\lambda;z) = 0, \qquad
    \widehat{\phi_i       [0]} V(\lambda;z) = \lambda(H_i) V(\lambda;z),
\label{V-highest-weight}
\end{equation}
for $\alpha\in\Delta_+$, $i=1,\dots,l$, $m\geqq0$, $n>0$.
\end{lem}
The second statement is due to the operator product expansions:
\begin{equation}
    \beta_\alpha(z)  V(\lambda;w) \sim 0, \qquad
    \gamma^\alpha(z) V(\lambda;w) \sim 0, \qquad
    \phi_i(z)        V(\lambda;w) \sim 
    \frac{\lambda(H_i)}{z-w} V(\lambda;w).
\label{ope:vo}
\end{equation}

Hence, the universality of Fock representations implies
\begin{cor}
\label{state-operator}
There exists a unique surjective homomorphism
\begin{equation}
    \Phi_\lambda: \Waki_{\lambda,k} \to \O_\lambda,
\label{state-operator-hom}
\end{equation}
which maps $\ghvac\tensor|\lambda\bosket$ to $V(\lambda;z)$. 
\end{cor}

It follows immediately from \eqref{def:operator-realization} that the
$\g$-submodule $\Waki^0_{\lambda,k}$ defined by
\eqref{def:wakimoto-finite} is isomorphically mapped to
\begin{equation}
    \calP_\lambda :=
    \{\, P(\gamma(z))V(\lambda;z)\,|\,
    P(x)\text{ is a polynomial of }x=(x^\alpha)_{\alpha\in\Delta_+}\,
    \}
\label{def:P-lambda}
\end{equation}
by $\Phi_\lambda$:
\begin{equation}
    \Phi_\lambda(P(\gamma[0])\ghvac\tensor|\lambda\bosket)
    =
    P(\gamma(z))V(\lambda;z).
\label{state-operator-hom:finite}
\end{equation}

\medskip
Since $\hat\g\oplus\Vir$ is realized in terms of ghosts and bosons through
\eqref{wakimoto-hom-km} and \eqref{wakimoto-hom-vir}, we can define action
of $X[m]$ ($X\in\g$, $m\in\Integer$) and $T[m]$ ($m\in\Integer$) on
$\O_\lambda$ by replacing $\beta_\alpha[n]$ etc.\ in $\omega(X[m])$,
$\omega(T[m])$ with $\widehat{\beta_{\alpha}[n]}$ etc.\ defined by
\eqref{def:operator-realization} respectively. In fact, their actions are 
described more simply, thanks to the following lemma: Assume that
$B_i(z)=\sum_{m\in\Integer}B_i[m]z^{-m-h_i}$ ($i=1,2$) are fields which
have the operator product expansion
\begin{equation}
    B_1(z) B_2(w)
    =
    \sum_{j=1}^N \frac{B_{12,j}(w)}{(z-w)^j}
    +
    \np B_1(z) B_2(w) \np,
\label{ope:B}
\end{equation}
where the normal order product $\np\ \np$ is defined by
\begin{equation}
    \np B_1[m] B_2[n] \np
    =
    \begin{cases}
    B_1[m] B_2[n], & m\leqq -h_1, \\
    B_2[n] B_1[m], & m>     -h_1
    \end{cases}
\label{def:normal-order:B}
\end{equation}
and the field $\np B_1(z) B_2(w)\np$ has no singularity at $z=w$. We
denote its restriction to the diagonal $\{z=w\}$ by $B_3(z)$:
\begin{equation}
    B_3(z) := \sum_{m\in\Integer} B_3[m] z^{-m-h_3}, \qquad
    B_3[m] = \sum_{n\in\Integer} \np B_1[n] B_2[m-n] \np,
\label{def:B3}
\end{equation}
where $h_3 = h_1 + h_2$.
\begin{lem}
For any $\Phi(z) \in \O_\lambda$, we have
$$
    \widehat{B_3[m]} \Phi(z) =
    \sum_{n\in\Integer} \NP\widehat{B_1[n]} \widehat{B_2[m-n]}\NP
    \Phi(z),
$$
where
$$
    \widehat{B_i[m]} \Phi(z) =
    \Res_{\zeta=z} (\zeta-z)^{m+h_i-1} B_i(\zeta) \Phi(z),
$$
and the normal ordering $\NP\ \NP$ is defined by the same rule as in
\eqref{def:normal-order:B}.
\end{lem}

\begin{cor}
For $X\in\g$, $m\in\Integer$ and $\Phi(z) \in \O_\lambda$, we have
\begin{align}
    \widehat{X[m]} \Phi(z) = 
    \Res_{\zeta=z} (\zeta-z)^m X(\zeta) \Phi(z),
\label{operator-realization:km}
\\
    \widehat{T[m]} \Phi(z) = 
    \Res_{\zeta=z} (\zeta-z)^{m+1} T(\zeta) \Phi(z),
\label{operator-realization:vir}
\end{align}
where action in the left hand side is defined by replacing
$\beta_\alpha[m]$ etc.\ in $\omega(X[m])$ and $\omega(T[m])$ by
$\widehat{\beta_\alpha[m]}$ etc.\ defined by
\eqref{def:operator-realization}. 
\end{cor}

Using the operator product expansions,
\begin{equation}
    T(z) V(\lambda;w) \sim
    \frac{\partial V(\lambda;w)}{z-w} + 
    \frac{\Delta_\lambda V(\lambda;w)}{(z-w)^2},
\label{ope:em-vo}
\end{equation}
where $\Delta_\lambda = (\lambda|\lambda+2\rho)/2\kappa$ is the conformal
weight of $V(\lambda;w)$, and the fact
\begin{equation}
    [T[-1], x(z)] = \partial x(z),
\label{com-rel:T[-1]}
\end{equation}
for $x=\beta_\alpha, \gamma^\alpha$ or $\partial\phi_i$, which is a direct
consequence of \eqref{def:T(z)}, we can prove the following formula by
induction.

\begin{lem}
For any $\Phi(z) \in \O_\lambda$, we have
\begin{equation}
    \partial \Phi(z) = [T[-1], \Phi(z)] = \widehat{T[-1]} \Phi(z).
\label{T[-1]-vo}
\end{equation}
\end{lem}

\subsection{Screening operators}
\label{screening}

Bosonization of the operator $\Scr(X;x,\partial_x)$ in
\eqref{nilp-left-action} gives the screening operator. The ghost sector of
the screening operator is defined for any positive root
$\alpha\in\Delta_+$ as follows:
\begin{equation}
    \Scr_\alpha(z) := \np \Scr(e_\alpha;\gamma(z), \beta(z)) \np.
\label{def:Scr(z)}
\end{equation}
The {\em screening operator} or the {\em screening current} is defined
for simple roots $\alpha_i$ ($i=1,\dots,l$) as the product of
$\Scr_{\alpha_i}$ and a bosonic vertex operator (see
\eqref{def:boson-vo}):
\begin{equation}
    \scr_i(z) := \Scr_{\alpha_i}(z) V(-\alpha_i;z)
    \in \O_{-\alpha_i}.
\label{def:scr(z)}
\end{equation}

Important property of screening currents is the following operator product
expansions. 
\begin{align}
    X(z) \scr_i(w) &\sim 0,
\label{ope:borel-scr}
\\
    F_j(z) \scr_i(w) &\sim
    - \kappa \delta_{i,j} 
      \frac{\partial}{\partial w} \frac{V(-\alpha_i;w)}{z-w},
\label{ope:F-scr}
\\
    T(z) \scr_i(w) &\sim 
    \frac{\partial}{\partial w} \frac{\scr_i(w)}{z-w},
\label{ope:em-scr}
\end{align}
where $X \in\borel_+$ and $i,j = 1,\dots,l$.  The ghost sector of the
screening current has the following operator product expansion, which
shall be used in computing explicit forms of the integral representations
of solutions of the KZB equations:
\begin{align}
    H^\gh(z) \Scr_\alpha(w)
    &\sim
    \frac{\alpha(H)\Scr_\alpha(w)}{z-w},
\label{ope:H-Scr}
\\
    T^\gh(z) \Scr_\alpha(w)
    &\sim
    \frac{\partial \Scr_\alpha(w)}{z-w} 
    +
    \frac{\Scr_\alpha(w)}{(z-w)^2},
\label{ope:em-Scr}
\\
    \Scr_\alpha(z) \Scr_\beta(w) 
    &\sim 
    \frac{f^{\alpha+\beta}_{\alpha,\beta} \Scr_{\alpha+\beta}(w)}{z-w},
\label{ope:Scr-Scr}
\\
    \Scr_\alpha(z) P(\gamma(w))
    &\sim
    \frac{(\Scr_\alpha P)(\gamma(w))}{z-w},
\label{ope:Scr-P}
\end{align}
for any $P(x)\in\Complex[x]$ ($x=(x^\alpha)_{\alpha\in\Delta}$), where
$f^{\alpha+\beta}_{\alpha,\beta}$ is the structure constant of the Lie
algebra $\n_+$, $[e_\alpha,e_\beta]=f_{\alpha,\beta}^{\alpha+\beta}
e_{\alpha+\beta}$, and
$(\Scr_\alpha P)(x)\in\Complex[x]$ is the polynomial $(\Scr_\alpha P)(x) =
\Scr(e_\alpha; x, \partial_x) P(x) \in\Complex[x]$.

\section{WZW models on elliptic curves}
\label{wzw}

In this section we recall the definition (or a characterization) of
$N$-point functions on elliptic curves.

\subsection{Space of conformal coinvariants and space of conformal blocks}
\label{cc-cb}

$N$-point functions of the WZW model take values in the space of the
conformal blocks which is the dual of the space of the conformal
coinvariants. (Exactly speaking, $N$-point functions are sections of a
vector bundle, a fiber of which is the space of conformal blocks. See
\secref{N-point-function}.) To define the space of conformal coinvariants
and conformal blocks, we first need a Lie algebra bundle over an elliptic
curve with marked points.

For each $q \in \Complex^\times$, $|q|<1$, we define an elliptic curve
$X = X_q$ by
\begin{equation}
    X_q := \Ctimes/q^\Integer,
\label{def:ell-curve}
\end{equation}
where $q^\Integer = \{ q^n \,|\,n\in\Integer \}$ is a multiplicative group
acting on $\Ctimes$ by $z \mapsto q^nz$. A Lie algebra bundle $\g^H$ is
defined for each $H\in\h$ by
\begin{equation}
   \g^H = \Ctimes \times \g/\sim,
\label{def:gH}
\end{equation}
where the equivalence relation $\sim$ is
\begin{equation}
    (t, A) \sim (qt, e^{- \ad H}A).
\label{def:eq-rel-gH}
\end{equation}

This Lie algebra bundle $\g^H$ has a natural connection, $\nabla_{d/dt} =
t d/dt$, and is decomposed into a direct sum of line bundles and a trivial
bundle with fiber $\h$:
\begin{equation}
  \g^H \simeqq \bigoplus_{\alpha\in\Delta} L_{\alpha(H)}
                 \oplus (\h \times X).
\label{decomp:gH}
\end{equation}
Here the line bundle $L_c$ ($c\in\Complex$) on $X$ is defined by
\begin{equation}
    L_c := (\C \times \C)/{\approx_c},
\label{def:L-c}
\end{equation}
where $\approx_c$ is an equivalence relation defined by
\begin{equation}
    (t, x) \approx_c (qt, e^{-c} x).
\label{def:approx-c}
\end{equation}

As usual, the structure sheaf on $X=X_q$ is denoted by $\O_X$ and the sheaf
of meromorphic functions on $X$ by $\K_X$. A stalk of a sheaf $\F$ on $X$
at a point $P \in X$ is denoted by $\F_P$.  When $\F$ is an $\O_X$-module,
we denote its fiber $\F_P/\m_P\F_P$ by $\F|_P$, where $\m_P$ is the maximal
ideal of the local ring $\O_{X,P}$.  Denote by $\F^\wedge_P$ the
$\m_P$-adic completion of $\F_P$.

We shall use the same symbol for a vector bundle and for a locally free
coherent $\O_X$-module consisting of its local holomorphic sections. For
instance, the invertible sheaf associated to the line bundle $L_c$ is also
denoted by the same symbol $L_c$. Denote by $\Omega_X^1$ the sheaf of
holomorphic 1-forms on $X$, which is isomorphic to $\O_X$ since $X$ is an
elliptic curve.  The fiberwise Lie algebra structure of the bundle $\g^H$
induces that of the associated sheaf $\g^H$ over $\O_X$.  Define the
invariant $\O_X$-inner product on $\g^H$ by
\begin{equation}
  (A|B) := \frac{1}{2 h^\vee} \Tr_{\g^H}(\ad A \ad B) \in \O_X
  \quad
  \text{for $A,B\in\g^H$},
  \label{def:inner-prod-gH}
\end{equation}
where the symbol $\ad$ denotes the adjoint representation of the $\O_X$-Lie
algebra $\g^H$. Then the inner product on $\g^H$ is invariant under the
translations with respect to the connection 
$\nabla:\g^H\to\g^H\tensor_{\O_X}\Omega_X^1$:
\begin{equation}
  d(A|B) = (\nabla A|B) + (A|\nabla B) \in \Omega_X^1
  \quad
  \text{for $A,B\in\g^H$}.
  \label{compati-conn-innerprod}
\end{equation}
Under the trivialization of $\g^H$ defined by the construction
\eqref{def:gH}, the connection $\nabla$ coincides with the exterior
derivative by $t\,d/dt$.

The fiber of $\g^H$ is isomorphic to $\g$. For any point $P$ on $X$
with $t(P)=z$, we put
\begin{equation}
    \g^P := (\g^H \tensor_{\O_X} \K_X)^{\wedge}_P,
\label{def:loop-alg}
\end{equation}
which is a topological Lie algebra non-canonically isomorphic to the loop
algebra $\g((t-z))$. Its subspace
$\g^P_+:=(\g^H)_P^{\wedge}\simeqq\g[[t-z]]$ is a maximal linearly compact
subalgebra of $\g^P$ under the $(t-z)$-adic linear topology.

Let us fix mutually distinct points $P_1, \ldots, P_N$ on $X$ whose
coordinates are $t = z_1, \ldots, z_N$ and put $D:=\{P_1,\dots,P_N\}$. We
shall also regard $D$ as a divisor on $X$ (i.e., $D=P_1+\cdots+P_N$). Denote
$X\setminus D$ by $\dot X$.  The Lie algebra
$\g^D:=\bigoplus_{i=1}^N\g^{P_i}$ has the natural 2-cocycle defined by
\begin{equation}
  \caff(A,B) := \sum_{i=1}^N \Res_{t = z_i} (\nabla A_i | B_i),
\label{def:cocycle}
\end{equation}
where $A=(A_i)_{i=1}^N, B=(B_i)_{i=1}^N \in \g^D$ and $\Res_{t=z}$ is the
residue at $t=z$. (The symbol ``$\caff$'' stands for ``Cocycle defining the
Affine Lie algebra''.)  We denote the central extension of $\g^D$ with respect
to this cocycle by $\hat\g^D$:
\begin{equation*}
    \hat\g^D := \g^D \oplus \C \khat,
\end{equation*}
where $\khat$ is a central element.  Explicitly the bracket of $\hat\g^D$ is
represented as
\begin{equation}
  [A, B] = ([A_i, B_i]^0)_{i=1}^N \oplus \caff(A,B) \khat
  \quad
  \text{for $A,B\in\g^D$,}
  \label{def:aff-alg-str}  
\end{equation}
where $[A_i, B_i]^0$ are the natural bracket in $\g^{P_i}$.  The Lie
algebra $\hat\g^P$ for a point $P$ is non-canonically isomorphic to the
affine Lie algebra $\hat\g$ and $\hat\g^{P_i}$ can be regarded as a
subalgebra of $\hat\g^D$.  Put $\g^P_+:=(\g^H)^\wedge_P$ as above.  Then
$\g^{P_i}_+$ can be also regarded as a subalgebra of $\hat\g^{P_i}$ and
$\hat\g^D$.

Let $\gHDpr$ be the space of global meromorphic sections of $\g^H$ which are
holomorphic on $\dot X$:
\begin{equation*}
    \gHDpr := \Gamma(X, \g^H(\ast D)).
\end{equation*}
There is a natural linear map from $\gHDpr$ into $\g^D$ which maps a
meromorphic section of $\g^H$ to its germ at $P_i$'s. The residue theorem
implies that this linear map is extended to a Lie algebra injective
homomorphism from $\gHDpr$ into $\hat\g^D$, which allows us to regard
$\gHDpr$ as a subalgebra of $\hat\g^D$.

\begin{defn}
\label{def:CC-CB}
The space of {\em conformal coinvariants} $\CC_H(X_q,D,M)$ and that of
{\em conformal blocks} $\CB_H(X_q,D,M)$ associated to
$\hat\g^{P_i}$-modules $M_i$ with the same level $\khat = k$ are
defined to be the space of coinvariants of $M := \bigotimes_{i=1}^N
M_i$ with respect to $\gHDpr$ and its dual:
\begin{equation}
    \CC_H(X_q,D,M) := M/\gHDpr M, \quad
    \CB_H(X_q,D,M) := (M/\gHDpr M)^\ast.
\label{def:conf-block}
\end{equation}
(In \cite{tu-ue-ya:89}, $\CC_H(X_q,D,M)$ and $\CB_H(X_q,D,M)$ are
called the space of {\em covacua} and that of {\em vacua}
respectively.)  In other words, the space of conformal coinvariants
$\CC_H(X_q,D,M)$ is generated by $M$ with relations
\begin{equation}
    A_\prin v \equiv 0
    \quad\text{for all $A_\prin \in \gHDpr$ and $v \in M$,}
\label{ward-cc}
\end{equation}
and a linear functional $\Phi$ on $M$ belongs to the space of conformal
blocks $\CB_H(X_q,D,M)$ if and only if it satisfies that
\begin{equation}
    \Phi(A_\prin v) = 0
    \quad\text{for all $A_\prin \in \gHDpr$ and $v \in M$.}
\label{ward-cb}
\end{equation}
These equations \eqref{ward-cc} and \eqref{ward-cb} are called  the
{\em Ward identities}.
\end{defn}

\subsection{$N$-point functions}
\label{N-point-function}

$N$-point functions are flat sections of a sheaf of conformal blocks
over the base space $S$ of a family $\Xtilde$ of pointed curves with
marked points defined as follows:
\begin{align*}
  & S :=
  \{\,(z;q;H)=(z_1,\ldots,z_N;q;H)
      \in(\Complex^\times)^N\times\Complex^\times\times\h \mid
      z_i/z_j \not\in q^\Z \ \text{if}\ i \ne j \,\},
  \\
  & \Xtilde := S\times\C^\times.
\end{align*}
Let $\tilde\pr=\pr_{\Xtilde/S}$ be the projection from $\Xtilde$ onto
$S$ along $\C^\times$, $\pi(z;q;H;t) = (z;q;H)$, and $\tilde p_i$ the
section of $\tilde\pr$ given by $z_i$:
\begin{equation*}
  \tilde p_i(z;q;H) := (z;q;H;z_i) \in \Xtilde
  \quad
  \text{for $(z;q;H)=(z_1,\ldots,z_N;q;H)\in S$}.
\end{equation*}

A {\em family of $N$-pointed elliptic curves} $\pr:\X\onto S$ is
constructed as follows. Define the action of $\Z$ on $\Xtilde$ by
\begin{equation}
  m\cdot(z;q;H;t) := (z;q;H;q^m t)
  \quad
  \text{for $m\in\Z$, $(z;q;H;t)\in\Xtilde$.}
  \label{def:action-Z-on-Xtilde}
\end{equation}
Let $\X$ be the quotient space of $\Xtilde$ by the action of $\Z$:
\begin{equation}
  \X := \Z\backslash\Xtilde.
  \label{def:X}
\end{equation}
Let $\pr_{\Xtilde/\X}$ be the natural projection from $\Xtilde$ onto $\X$ and
$\pr=\pr_{\X/S}$ the projection from $\X$ onto $S$ induced by $\tilde\pr$. 
We put
\begin{equation*}
  p_i:=\pr_{\Xtilde/\X}\circ\tilde p_i, \quad
  P_i := p_i(S), \quad
  D := \bigcup_{i=1}^N P_i, \quad
  \dot\X := \X \setminus D, \quad
  \tilde D := \pr_{\Xtilde/\X}^{-1}(D).
\end{equation*}
Here $p_i$ is the section of $\pr$ induced by $\tilde p_i$ and $D$ is also
regarded as a divisor $\sum_{i=1}^N P_i$ on $\X$.  
The fiber of $\pr$ at $(z;q;H)=(z_1,\ldots,z_N;q;H)\in S$
is an elliptic curve with modulus $q$ and marked points $z_1,\ldots,z_N$.

We refer to \cite{fel-wie:96} for the construction of a sheaf of conformal 
blocks and its flat connection and, using their result, define the
$N$-point functions as follows. (See also \cite{suz:96}.)

We identify each fiber of $\g^H$ with $\g$ via the standard trivialization
defined by the construction of $\g^H$, \eqref{def:gH}. Then the algebra
$\g^P$ defined by \eqref{def:loop-alg} is identified with $\g((t-z))$
where $t$ is the coordinate on the complex plane and $z$ is the coordinate
of $P$. For $X\in\g$, the element $X\tensor (t-z)^m$ of $\g^P \simeqq
\g((t-z))$ is denoted by $X[m]$. The Virasoro generator defined by the
Sugawara construction \eqref{em-tensor} is denoted by $T[m]$:
\begin{equation}
    T[m]
    =
    \frac{1}{2\kappa} 
    \sum_{p=1}^{\dim \g} \sum_{n\in\Integer} \NP J_p[m-n] J^p[n]\NP.
\label{sugawara}
\end{equation}

Let us denote the representation of $\hat\g^{P_i}$ on the $i$-th
component of the tensor product $M=\bigotimes_{i=1}^N M_i$ by $\rho_i$
and its dual by $\rho^\ast_i$: for $v\in M$, $v^\ast\in M^\ast$,
$A\in\hat\g^{P_i}$,
\begin{equation}
    \langle \rho_i^\ast(A) v^\ast, v \rangle
    =
    - \langle v^\ast, \rho_i(A) v \rangle
\label{def:dual-rep}
\end{equation}
where $\langle\cdot, \cdot\rangle$ is the pairing of $M^\ast$ and $M$. 
We assume that the Virasoro algebra with central charge $\cvir$ acts on
$M_i$ and $M_i^\ast$ through the Sugawara construction \eqref{sugawara}.

For $v\in M$, $v^\ast\in M^\ast$, $H\in\h$ and a multi-valued meromorphic
function $f(t)$ on $X_q$, all poles of which belong to $\{z_1, \dots,
z_N\} \pmod {q^\Integer}$, we define
\begin{alignat}{2}
    \rho\left( T\left\{ f(t) \frac{d}{dt}\right\}\right) v &:=
    \sum_{i=1}^N \sum_{m\in\Integer} f_{i,m+1} \rho_i(T[m]) v,\qquad &
    \rho^\ast\left( T\left\{ f(t) \frac{d}{dt}\right\}\right) v^\ast &:=
    \sum_{i=1}^N \sum_{m\in\Integer} f_{i,m+1} \rho^\ast_i(T[m]) v^\ast,
\label{def:T(vec-field)}
\\
    \rho( H\{ f(t) \}) v &:=
    \sum_{i=1}^N \sum_{m\in\Integer} f_{i,m} \rho_i(H[m]) v,\qquad &
    \rho^\ast( H\{ f(t) \}) v^\ast &:=
    \sum_{i=1}^N \sum_{m\in\Integer} f_{i,m} \rho^\ast_i(H[m]) v^\ast,
\label{def:H(function)}
\end{alignat}
where $f(t) = \sum_{m\in\Integer} f_{i,m}(t-z_i)^m$ is the Laurent
expansion of $f(t)$ around $t=z_i$.

Fix a meromorphic function $Z(z;q;H;t)$ on $\Xtilde$ with poles only at
$\{z_1, \dots,z_N\} \pmod {q^\Integer}$ (namely
$Z(z;q;H;t)\in\Gamma(\Xtilde, \O_\Xtilde(\ast\tilde D))$) satisfying
\begin{equation}
    Z(z;q;H;qt) = Z(z;q;H;t) - 1.
\label{def:Z}
\end{equation}
We abbreviate $Z(q;z;H;t)$ as $Z(t)$.

\begin{example}
We may take the following function as $Z(t)=Z(z;q;H;t)$:
\begin{equation}
    Z(z;q;H;t) 
    := \frac{1}{2\pi i} \frac{t}{z_{i_0}}
       \frac{d}{dt} \log\theta_{11}(t/z_{i_0};q),
\end{equation}
for $i_0\in\{1,\dots,N\}$.
\end{example}

Let us take a coordinate system of $\h$ as $\h\owns H = \sum_{a=1}^l \xi_a
H_a$, where $\{H_a\}_{a=1}^l$ is an orthonormal basis of $\h$.

\begin{defn}
\label{def:N-point-function}
A multi-valued holomorphic function $\Psi(z;q;H)$ on $S$ with values in
$M^\ast$ is called an {\em $N$-point function} in genus one if it
satisfies the following conditions (I), (II), (III) and (IV).:

\begin{enumerate}
\renewcommand{\labelenumi}{(\Roman{enumi})}
\item
For any $(z,q,H) \in S$, $\Psi(z;q;H) \in \CB_H(X_q,D,M)$;

\item
For $j=1, \dots, N$,
\begin{equation}
    \frac{\partial}{\partial z_j} \Psi(z;q;H) =
    \rho^\ast_j(T[-1]) \Psi(z;q;H);
\label{N-point:eq1}
\end{equation}

\item
\begin{equation}
    \left(q \frac{\partial}{\partial q} + \frac{\cvir}{24} \right)
    \Psi(z;q;H) =
    \rho^\ast\left(T\left\{Z(t)t\frac{d}{dt}\right\}\right) \Psi(z;q;H);
\label{N-point:eq2}
\end{equation}

\item
For $r=1, \dots, l$,
\begin{equation}
    \frac{\partial}{\partial \xi_r} \Psi(z;q;H) =
    - \rho^\ast(H_r\{Z(t)\}) \Psi(z;q;H).
\label{N-point:eq3}
\end{equation}
\end{enumerate}
\end{defn}

Assume that each $M_i$ contains a $\g$-submodule $V_i$
such that for any $m>0$ and $X\in\g$,
\begin{equation}
    X[m] V_i = 0,
\label{nilp=zero}
\end{equation}
and the Casimir operator $C_2=\sum_p J_p J^p$ acts as a multiplication,
\begin{equation}
    \rho_i(C_2)
    =
    c_2^{(i)} \id_{V_i}.
\label{casimir-on-Vi}
\end{equation}
Let us define a function $\Pi(q;H)$ by
\begin{equation}
    \Pi(q;H) :=
    q^{\dim\g/24} (q;q)_\infty^l
    \prod_{\alpha\in\Delta_+} 2\sinh(\alpha(H)/2)
    \prod_{\alpha\in\Delta} (q e^{\alpha(H)};q)_\infty.
\label{def:Pi}
\end{equation}
Then we can restrict an $N$-point function $\Psi(z;q;H)$ to
$V=\bigotimes_{i=1}^N V_i$:
\begin{equation}
    \Psi_V(z;q;H) := \Psi(z;q;H)|_V,
\label{def:Psi-V}
\end{equation}
and $\tilde\Psi_V(z;q;H)=\Pi(q;H)\Psi_V(z;q;H)$ satisfies the following
system, which is called the {\em Knizhnik-Zamolodchikov-Bernard} (KZB)
equations first found by Bernard \cite{ber:88-1}. (See Theorem~4.1 of
\cite{fel-wie:96}.): The functions $\sigma_c(z)$ and $\zeta(z)$ below are
defined by \eqref{def:sigma} and \eqref{def:zeta}.
\begin{enumerate}
\renewcommand{\labelenumi}{(\Roman{enumi}${}'$)}
\item
For any $H\in\h$,
\begin{equation}
    \sum_{i=1}^N \rho_i^\ast(H) \tilde\Psi_V(z;q;H) = 0;
\label{kzb:weight0}
\end{equation}

\item
For $j=1,\dots,N$, 
\begin{equation}
    \kappa\left(
    z_j \frac{\partial}{\partial z_j} + \frac{c_2^{(j)}}{2\kappa}
    \right)
    \tilde\Psi_V(z;q;H)
    =
    \sum_{r=1}^l
    \rho_j^\ast(H_r)
    \frac{\partial}{\partial \xi_r} \tilde\Psi_V(z;q;H)
    +
    \sum_{i\neq j}
    (\rho_i^\ast \tensor \rho_j^\ast \Omega(z_i,z_j))
    \tilde\Psi_V(z;q;H),
\label{kzb}
\end{equation}
where
\begin{equation}
    \Omega(z,w)
    =
    \Omega(z,w;q;H)
    :=
    - \sum_{\alpha\in\Delta}
    \sigma_{-\alpha(H)}\left(\frac{z}{w}\right)
    e_\alpha \tensor e_{-\alpha}
    - \sum_{r=1}^l
    \zeta\left(\frac{z}{w}\right) H_r \tensor H_r.
\label{def:Omega}
\end{equation}

\item
\begin{equation}
    2\kappa q\frac{\partial}{\partial q} \tilde\Psi_V(z;q;H)
    =
    \sum_{r=1}^l
    \left(\frac{\partial}{\partial\xi_r}\right)^2 \tilde\Psi_V(z;q;H)
    +
    \sum_{i,j=1}^N
    (\rho_i^\ast\tensor \rho_j^\ast H(z_i,z_j))
    \tilde\Psi_V(z;q;H),
\label{kzb-heat}
\end{equation}
where
\begin{multline}
    H(z,w)
    =
    H(z,w;q;H)
    :=\\:=
    - \sum_{\alpha\in\Delta}
    \left(
     \zeta\left(\frac{e^{\alpha(H)}z}{w}\right)
     -
     \zeta(e^{\alpha(H)})
    \right)
    \sigma_{-\alpha(H)}\left(\frac{z}{w}\right)
    e_\alpha \tensor e_{-\alpha}
    - \sum_{r=1}^l
    \frac12 
    \left(
     \zeta\left(\frac{z}{w}\right)^2
     +
     \frac{z}{w} \zeta'\left(\frac{z}{w}\right)
    \right) 
    H_r \tensor H_r.
\label{def:H}
\end{multline}
\end{enumerate}

Conversely, restriction of an $N$-point function to
$V$ is characterized by the KZB equations with additional conditions. For
example, Felder and Wieczerkowski \cite{fel-wie:96} used the automorphic
properties and asymptotic behavior of $\tilde\Psi_V$ as the additional
conditions, while Suzuki \cite{suz:96} found a holonomic system
characterizing $\tilde\Psi_V$ which includes the KZB equations.

\section{$N$-point functions from Wakimoto modules}
\label{N-point-function:wakimoto}

In this section we construct $N$-point functions for Wakimoto modules,
$M_i=\Waki_{\lambda_i,k}$, where $\lambda_i\in\h^\ast$ and
$\sum_{i=1}^N\lambda_i$ belongs to the positive root lattice of $\g$. 

Fix ordered sets of the simple roots of $\g$, $\{\alpha_{i_1}, \ldots,
\alpha_{i_M}\}$, such that
\begin{equation}
    \sum_{j=1}^M \alpha_{i_j} = \sum_{i=1}^N \lambda_i.
\label{charge-conservation}
\end{equation}

Denote the following linear map by $\psi(z;t;q;H)$, where
$z=(z_1,\dots,z_N)$, $t=(t_1,\dots,t_M)$ and $(z,t)$ belongs to (the
universal covering of)
$(\Complex^\times)^{N+M}\setminus\{\,$diagonals$\,\}$:
\begin{multline}
    \Waki_{\lambda_1,k}\tensor \cdots \tensor \Waki_{\lambda_N,k}
    \owns v_1 \tensor \cdots \tensor v_N \mapsto
\\
    \mapsto
    \Tr_{\Waki_{\mu,k}}(
    \Phi_{\lambda_1}(v_1;z_1) \cdots \Phi_{\lambda_N}(v_N;z_N)
    \scr_{i_1}(t_1) \cdots \scr_{i_M}(t_M)
    q^{T[0]-\cvir/24} e^{H[0]}
    ) dt_1 \wedge \cdots dt_M.
\label{def:psi}
\end{multline}
Note that, thanks to \eqref{O-lambda-act-on-wakimoto} and
\eqref{charge-conservation}, the operator inside the bracket in the right
hand side is an endomorphism of $\Waki_{\mu,k}$.

\begin{prop}
\label{local-system-def-by-psi}
There is a local system of rank one $\Line$ on $\{\,(z;t;q;H)\,|\,
(z;t)=(z_1,\dots,z_N;t_1,\dots,t_M)\in X_q^{N+M} \setminus
\{\text{diagonals}\}, q\in \Complex^\times, H\in \h\,\}$ such that
$\psi(z;t;q;H)(v)$ is a holomorphic section of $\Line$ for any $v\in
M$. Namely, the monodromy of $\psi(z;t;q;H)(v)$ is independent of $v$.
\end{prop}

\begin{proof}
First note that $\Phi_{\lambda_i}(v_i;z_i)\in\O_\lambda$ is of the form,
\begin{equation}
\begin{split}
    \Phi_{\lambda_i}(v_i;z_i) &=
    \sum
    \widehat{x_1[m_1]} \cdots \widehat{x_n[m_n]} V(\lambda;z)
\\
    &=
    \sum
    \oint_{C_{i1}} d\zeta_{i1} (\zeta_{i1}-z_i)^{m_1+h_1-1}
    \cdots
    \oint_{C_{in}} d\zeta_{in} (\zeta_{in}-z_i)^{m_n+h_n-1}
    x_1(\zeta_{i1}) \cdots x_n(\zeta_{in}) V(\lambda_i;z_i),
\end{split}
\label{Phi-lambda:form}
\end{equation}
where $x_j$ is $\beta_\alpha$, $\gamma^\alpha$ or $\partial\phi_j$ and
$m_j\in\Integer$. The contour $C_{ij}$ encircles $z_i$, lies outside of
$C_{ij'}$ ($j'>j$) and does not contain $0$ and $z_{i'}$ ($i'\neq i$)
inside it.

It follows from \eqref{def:scr(z)}, \eqref{def:Scr(z)} and this expression 
\eqref{Phi-lambda:form} that $\psi(z;t;q;H)(v_1\tensor \cdots \tensor
v_N)$ is sum of integrals of the form
\begin{gather}
    \oint_{C_{i1}} d\zeta_{i1} \cdots \oint_{C_{in}} d\zeta_{in}
    ( \text{rational function of } \zeta_{ij}, z_i)
    F^\gh(\zeta_{ij}, t_i; q; H) F^\bos(\zeta_{ij}, z_i, t_i; q; H)
    dt_1 \wedge \cdots dt_M,
\label{psi:parts}
\\
    F^\gh(\zeta_{ij}, t_i; q; H)
    :=
    \Tr_{\ghFock}\bigl(
        (\text{polynomial of }
        \beta_\alpha(\zeta_{ij}), \beta_\alpha(t_i),
        \gamma^\alpha(\zeta_{ij}), \gamma^\alpha(t_i) )
        q^{T^\gh[0]} e^{H^\gh[0]} 
        \bigr)
\label{def:Fgh}
\\
\begin{aligned}
    F^\bos(\zeta_{ij}, z_i, t_i; q; H)
    :=&
    \Tr_{\bosFock_\mu}\bigl(
        \prod_j \phi_{r_{1j}}(\zeta_{1j}) V(\lambda_1;z_1)
        \cdots
        \prod_j \phi_{r_{Nj}}(\zeta_{Nj}) V(\lambda_N;z_N) \times
    \\
    &\times 
        V(-\alpha_{i_1};t_1) \cdots V(-\alpha_{i_M};t_M)
        q^{T^\phi[0]} e^{\phi[H;0]}
        \bigr)
\end{aligned}
\label{def:Fbos}
\end{gather}
where ``polynomial of $\beta_\alpha$, etc.'' contain possibly normal
ordered products coming from $\Scr_{\alpha_{i_j}}(t_j)$.

The first trace in the integrand in \eqref{psi:parts},
$F^\gh(\zeta_{ij},t_i;q;H)$, has singularities at $\zeta_{ij}=\zeta_{i'j'}$
and $\zeta_{ij}=t_{i'}$, which are poles because of the operator product
expansions \eqref{ope:ghost}.

The second trace in the integrand in \eqref{psi:parts},
$F^\bos(\zeta_{ij},z_i,t_i;q;H)$, has singularities at (1) $\zeta_{ij} =
\zeta_{i'j'}$; (2) $\zeta_{ij} = z_{i'}$ or $t_{i'}$; (3) $z_i = z_{i'}$,
$z_i=t_{i'}$, $t_i=t_{i'}$. The first singularities (1) are poles because
of the operator product expansion \eqref{ope:boson}. The second
singularities (2) are also rational by virtue of the third expansion in
\eqref{ope:vo}. The formula (cf.\ \cite{kur:91} (5.32))
\begin{equation}
    V(\lambda;z) V(\lambda';w) =
    \np V(\lambda;z) V(\lambda';w) \np (z-w)^{(\lambda|\lambda')/\kappa}
\label{boson-vo:normal-order}
\end{equation}
implies that $F^\bos$ has non-trivial monodromy around the singularities
(3):
\begin{equation}
\begin{aligned}
    F^\bos(\zeta;z;t;q;H) 
    &\to e^{2\pi i(\lambda_i|\lambda_j)/\kappa} 
         F^\bos(\zeta;z;t;q;H),
    \text{ when $z_i$ goes around $z_j$, }(1\leqq i<j\leqq N),
\\
    F^\bos(\zeta;z;t;q;H) 
    &\to e^{2\pi i(\lambda_i|-\alpha_{i_j})/\kappa} 
         F^\bos(\zeta;z;t;q;H), 
    \text{ when $z_i$ goes around $t_j$, }
    (1\leqq i\leqq N, 1\leqq j\leqq M)),
\\
    F^\bos(\zeta;z;t;q;H) 
    &\to e^{2\pi i(\alpha_{i_j}|\alpha_{i_{j'}})/\kappa}
         F^\bos(\zeta;z;t;q;H), 
    \text{ when $t_{i_j}$ goes around $t_{i_{j'}}$, }
    (1\leqq j<j'\leqq M).
\end{aligned}
\label{monodromy:Fbos}
\end{equation}

Summarizing, we conclude that the integrand in \eqref{psi:parts} is
rational function with respect to $\zeta_{ij}$'s and also rational with
respect to $z_i$'s and $t_i$'s except at $z_i = z_{i'}$, $z_i=t_{i'}$,
$t_i=t_{i'}$, where it has the same monodromies as $F^\bos$,
\eqref{monodromy:Fbos}. 

As a next step, we show that $\psi(z;t;q;H)(v_1\tensor \cdots \tensor
v_N)$ has the same monodromy as the integrand in \eqref{psi:parts}. This
is proved by applying the following lemma iteratively.
\begin{lem}
\label{monodromy}
Assume that a function $F(\zeta, z, t)$ is rational with respect to
$\zeta$ and has monodromy with respect to $z$ and $t$ around the diagonal
$z=t$: $F(\zeta,z,t) \to c F(\zeta,z,t)$ when $z$ goes around $t$, where
$c$ is a constant. Then
\begin{equation}
    \psi(z,t) = \oint_{C(z)} F(\zeta,z,t) d\zeta
\label{def:psi-example}
\end{equation}
has the same monodromy around $z=t$, where $C(z)$ is a small contour
surrounding $z$.
\end{lem}

\begin{proof}
Fix a small circle $\gamma$ around $t$, 
$$
    \gamma(\theta) = t+\eps\exp(i\theta), \qquad
    \theta\in[0,2\pi].
$$
When $z$ goes around $t$ along $\gamma$, $F(\zeta,z,t)$ is multiplied by
$c$. Since $F$ is rational with respect to $\zeta$, the integration
contour $C(z)$ in \eqref{def:psi-example} can be replaced with a cycle
$\gamma_+ - \gamma_-$, where
$$
    \gamma_\pm(\theta) = t + \eps_\pm \exp(i\theta), \qquad
    \theta\in[0,2\pi]
$$
and $\eps_\pm$ are suitable constants satisfying $\eps_-<\eps<\eps_+$. Now 
that $\gamma_\pm$ do not depend on $z$ and do not intersect with $\gamma$, 
it is obvious that 
$$
    \psi(z,t) = \oint_{\gamma_+}F(\zeta,z,t) d\zeta
              - \oint_{\gamma_-}F(\zeta,z,t) d\zeta
$$
is multiplied by $c$ when $z$ goes around $t$ along $\gamma$.
\end{proof}

We can similarly prove that the monodromies of $\psi(z;t;q;H)(v_1\tensor
\cdots \tensor v_N)$ around the cycles of $X_q$ (along the paths,
$z_i=r\exp(2\pi i \theta)$ ($0\leqq \theta \leqq 2\pi$) and $z_i \to
qz_i$, and the same for $t_i$) do not depend on $v_1\tensor\cdots\tensor
v_N$. This completes the proof of the proposition.
\end{proof}

\begin{rem}
We can also write down an explicit expression of the second trace in the
integrand in \eqref{psi:parts}, $\Tr_{\bosFock_\mu}\bigl((\prod_{i,j}
\partial\phi_{r_{ij}}(\zeta_{ij}) V(\lambda_i;z_i) ) q^{T^\phi[0]}
e^{\phi[H;0]} \bigr)$, by using \eqref{boson-vo-1-loop:0-mode} and
\eqref{boson-vo-1-loop:non-0-mode}. In fact, since
$$
    \frac{\partial}{\partial z}
    \left. \frac{\partial}{\partial t}\right|_{t=0}
    \tilde V(t \lambda;z) 
    =
    \partial\phi(\lambda;z),
$$
applying differential operators of the form $\partial_z \partial_t|_{t=0}$
to \eqref{boson-vo-1-loop:non-0-mode} and combining the result with
\eqref{boson-vo-1-loop:0-mode}, we have a desired expression, which also
shows that this trace is rational with respect to $\zeta_{ij}$ and has
monodromy around the diagonals $z_i=z_{i'}$ etc.
\end{rem}

\begin{thm}
\label{int-rep-of-N-point-function}
\begin{equation}
    \Psi(z;q;H) = \int_{{\mathcal C}(z,q,H)} \psi(z;t;q;H)
\label{def:Psi}
\end{equation}
is an $N$-point function with values in
$(\Waki_{\lambda_1,k}\tensor\cdots\tensor\Waki_{\lambda_N,k})^\ast$, where
${\mathcal C}(z,q,H)$ is a family of $M$-cycles with coefficients in the
local system $\Line^\ast$ dual to $\Line$.
\end{thm}

\begin{rem}
We refer to \cite{aom-kit:94} or \cite{fel-var:95} for integrals over
cycles with coefficients in the local
system. \propref{local-system-def-by-psi} guarantees that the integration
in the right hand side of \eqref{def:Psi} is well-defined and that the
right hand sides of \eqref{N-point:eq1}, \eqref{N-point:eq2} and
\eqref{N-point:eq3} are meaningful.

\end{rem}

\begin{proof}
This can be shown in almost the same way as Proposition 3.4.1 of
\cite{suz:96}.

The condition (I) of \defref{def:N-point-function} is checked as
follows. Fix $(z,q,H)\in S$. The condition (I) means that for any
$J(t)\in\gHDpr$,
\begin{equation}
    \sum_{i=1}^N \rho_i^\ast(J(t)) \Psi(z;q;H) = 0.
\label{Psi-is-CB}
\end{equation}
Thanks to the decomposition \eqref{decomp:gH}, we may assume $J(t) = X
\tensor f(t)$ where $X\in\g_\alpha$ ($\alpha\in\Delta\sqcup\{0\}$,
$\g_0:=\h$), $f(t)\in \Gamma(X_q,L_{\alpha(H)}(\ast D))$ ($L_0=\O_{X_q}$). 
The left hand side of \eqref{Psi-is-CB} is equal to
\begin{equation}
\begin{split}
    &
    \int_{{\mathcal C}(z,q)} 
    \sum_{i=1}^N \rho_i^\ast(X\tensor f(t)) \psi(z;t;q;H)
    (v_1\tensor\cdots\tensor v_N)
\\
    =&
    \int_{{\mathcal C}(z,q)} 
    -\sum_{i=1}^N
    \Tr
    \bigl(
    \Phi_{\lambda_1}(v_1;z_1) \cdots 
    (X\tensor f(t))_{t=z_i}^\wedge \Phi_{\lambda_i}(v_i;z_i) \cdots
    \scr_{i_j}(t_j) \cdots q^{T[0]-\cvir/24} e^{H[0]}
    \bigr)
\\
    =&
    \int_{{\mathcal C}(z,q)} 
    -\sum_{i=1}^N
    \Res_{\zeta=z_i} f(\zeta)
    \Tr
    \bigl(
    \Phi_{\lambda_1}(v_1;z_1) \cdots 
    X(\zeta) \Phi_{\lambda_i}(v_i;z_i) \cdots \scr_{i_j}(t_j) \cdots
    q^{T[0]-\cvir/24} e^{H[0]}
    \bigr) d\zeta,
\end{split}
\label{Psi-is-CB:tmp1}
\end{equation}
where $\Tr$ is $\Tr_{\Waki_{\mu,k}}$ and $(X\tensor f(t))_{t=z_j}$ is the
Laurent expansion of $X\tensor f(t)$ at $t=z_j$. The last line is due to
the following fact, which is easily checked by
\eqref{operator-realization:km}: for any $\Phi(z)\in\O_\lambda$ and
$X\in\g$,
\begin{equation}
    (X\tensor f(t))_{t=z}^\wedge \Phi(z)
    =
    \Res_{\zeta=z} f(\zeta)X(\zeta) \Phi(z) d\zeta.
\label{action-of-Laurent-exp}
\end{equation}
By the commutativity of current $X(\zeta)$, vertex operator
$\Phi_\lambda(v;z)$ and screening current $\scr_i(z)$, we have
\begin{equation}
\begin{split}
    &
    f(\zeta)
    \Tr
    \bigl(
    \Phi_{\lambda_1}(v_1;z_1) \cdots 
    X(\zeta) \Phi_{\lambda_i}(v_i;z_i) \cdots \scr_{i_j}(t_j) \cdots
    q^{T[0]-\cvir/24} e^{H[0]}
    \bigr) d\zeta
\\
    =&
    f(\zeta)
    \Tr
    \bigl(X(\zeta) 
    \Phi_{\lambda_1}(v_1;z_1) \cdots 
    \Phi_{\lambda_i}(v_i;z_i) \cdots \scr_{i_j}(t_j) \cdots
    q^{T[0]-\cvir/24} e^{H[0]}
    \bigr) d\zeta
\\
    =&
    f(\zeta)
    \Tr
    \bigl(
    \Phi_{\lambda_1}(v_1;z_1) \cdots 
    \Phi_{\lambda_i}(v_i;z_i) \cdots \scr_{i_j}(t_j) \cdots
    q^{T[0]-\cvir/24} e^{H[0]}X(\zeta)
    \bigr) d\zeta
\\
    =&
    f(q\zeta)
    \Tr
    \bigl(
    \Phi_{\lambda_1}(v_1;z_1) \cdots 
    \Phi_{\lambda_i}(v_i;z_i) \cdots \scr_{i_j}(t_j) \cdots
    X(q\zeta) q^{T[0]-\cvir/24} e^{H[0]}
    \bigr) q d\zeta,
\end{split}
\label{Psi-is-CB:tmp2}
\end{equation}
where we used
\begin{equation}
    e^{H[0]} X(\zeta) = e^{\alpha(H)} X(\zeta) e^{H[0]}, \qquad
    q^{T[0]} X(\zeta) = q X(q\zeta) q^{T[0]},
\label{exp(ad(T[0])/exp(ad(H[0]))X(z)}
\end{equation}
and $f(qt) = e^{-\alpha(H)}f(t)$ (cf.\ \eqref{def:L-c}). Therefore,
$$
    f(\zeta) \Tr
    \bigl(\Phi_{\lambda_1}(v_1;z_1) \cdots 
    X(\zeta) \Phi_{\lambda_i}(v_i;z_i) \cdots \scr_{i_j}(t_j) \cdots
    q^{T[0]-\cvir/24} e^{H[0]} 
    \bigr) d\zeta
    \in
    \Gamma(X_q,\Omega^1_X(\ast D)).
$$
Hence the sum of its residues at $\zeta=z_j$ ($j=1,\dots,N$) is zero by
the residue theorem. Thus \eqref{Psi-is-CB:tmp1} implies
\eqref{Psi-is-CB}.

The condition (II) of \defref{def:N-point-function} is a direct
consequence of \eqref{T[-1]-vo}.

We prove the condition (III), assuming
$|q|<|t_M|<\cdots<|t_1|<|z_N|<\cdots<|z_1|<1$. The general case follows
from this case by the analytic continuation. By the same argument as
\eqref{Psi-is-CB:tmp1}, we have
\begin{multline}
    \rho^\ast\left(T\left\{Z(t)t\frac{d}{dt}\right\}\right) \psi(z;q;H)
    = \\ =
    \sum_{i=1}^N
    \Res_{\zeta=z_i} Z(\zeta) 
    \Tr
    \bigl(T(\zeta)
    \cdots \Phi_{\lambda_i}(v_i;z_i) \cdots
    \scr_{i_j}(t_j) \cdots q^{T[0]-\cvir/24} e^{H[0]}
    \bigr) d^M t \tensor \zeta d\zeta.
\label{T(Z(t))Psi:tmp1}
\end{multline}
Deforming the integration contour, the right hand side of
\eqref{T(Z(t))Psi:tmp1} is rewritten as
\begin{multline}
    \frac{1}{2\pi i}
    \left(
    \oint_{|\zeta| = 1} - \oint_{|\zeta|=|q|}
     - \sum_{j=1}^M \oint_{\zeta=t_j}
    \right)
    \zeta d\zeta Z(\zeta) \times 
    \\ \times
    \Tr
    \bigl(T(\zeta)
    \cdots \Phi_{\lambda_i}(v_i;z_i) \cdots
    \scr_{i_j}(t_j) \cdots q^{T[0]-\cvir/24} e^{H[0]}
    \bigr) d^M t
\label{T(Z(t))Psi:tmp2}
\end{multline}
The first integral in \eqref{T(Z(t))Psi:tmp2} is equal to 
\begin{equation}
\begin{split}
    &\frac{1}{2\pi i}
    \oint_{|\zeta| = 1} \!\!\!
    \zeta \, d\zeta \, Z(\zeta) 
    \Tr
    \bigl(
    \cdots \Phi_{\lambda_i}(v_i;z_i) \cdots
    \scr_{i_j}(t_j) \cdots q^{T[0]-\cvir/24} e^{H[0]} T(\zeta)
    \bigr) d^M t
\\
    =&\frac{1}{2\pi i}
    \oint_{|\zeta| = 1} \!\!\!
    \zeta \, d\zeta \, Z(\zeta) 
    \Tr
    \bigl(
    \cdots \Phi_{\lambda_i}(v_i;z_i) \cdots
    \scr_{i_j}(t_j) \cdots q^2 T(q\zeta) q^{T[0]-\cvir/24} e^{H[0]} 
    \bigr) d^M t
\\
    =&\frac{1}{2\pi i}
    \oint_{|\zeta| = 1} \!\!\!
    q\zeta \, d(q\zeta) \, Z(\zeta) 
    \Tr
    \bigl(
    \cdots \Phi_{\lambda_i}(v_i;z_i) \cdots
    \scr_{i_j}(t_j) \cdots T(q\zeta) q^{T[0]-\cvir/24} e^{H[0]} 
    \bigr) d^M t
\\
    =&
    \frac{1}{2\pi i}
    \oint_{|\zeta| = |q|} \!\!\!
    \zeta d\zeta \, Z(q^{-1}\zeta) 
    \Tr
    \bigl(
    \cdots \Phi_{\lambda_i}(v_i;z_i) \cdots
    \scr_{i_j}(t_j) \cdots T(\zeta) q^{T[0]-\cvir/24} e^{H[0]} 
    \bigr) d^M t,
\end{split}
\label{T(Z(t))Psi:tmp3}
\end{equation}
where we used
\begin{equation}
    e^{H[0]} T(\zeta) = T(\zeta) e^{H[0]}, \qquad
    q^{T[0]} T(\zeta) = q^2 T(q\zeta) q^{T[0]}.
\label{exp(ad(T[0])/exp(ad(H[0]))T(z)}
\end{equation}
The second integral in \eqref{T(Z(t))Psi:tmp2} is equal to
\begin{equation}
    \frac{1}{2\pi i}
    \oint_{|\zeta| = |q|}
    \zeta d\zeta \, Z(\zeta) 
    \Tr
    \bigl(
    \cdots \Phi_{\lambda_i}(v_i;z_i) \cdots
    \scr_{i_j}(t_j) \cdots T(\zeta) q^{T[0]-\cvir/24} e^{H[0]} 
    \bigr) d^M t
\label{T(Z(t))Psi:tmp4}
\end{equation}
The third integral in \eqref{T(Z(t))Psi:tmp2} turns into a term of the
form $\partial/\partial t_j (\cdots)$ by the operator product expansion
\eqref{ope:em-scr}, and therefore the sum of those terms is an exact
$M$-form. Hence they do not contribute to the integral over ${\mathcal
C}$. Thus by summing up \eqref{T(Z(t))Psi:tmp1}, \eqref{T(Z(t))Psi:tmp3}
and \eqref{T(Z(t))Psi:tmp4} and using the property of $Z$, \eqref{def:Z},
we obtain
\begin{equation}
\begin{split}
    \rho^\ast\left(T\left\{Z(t)t\frac{d}{dt}\right\}\right) \Psi(z;q;H)
    &=
    \int_{{\mathcal C}(z,q)}
    \Tr
    \bigl(
    \cdots \Phi_{\lambda_i}(v_i;z_i) \cdots
    \scr_{i_j}(t_j) \cdots T[0] q^{T[0]-\cvir/24} e^{H[0]} 
    \bigr)d^M t
\\
    &=
    \left(q \frac{\partial}{\partial q} + \frac{\cvir}{24}\right)
    \Psi(z;q;H),
\end{split}
\end{equation}
which proves the condition (III).

The condition (IV) is proved in the same way.
\end{proof}

Recall that the Wakimoto module $\Waki_{\lambda,k}$ contains a
$\g$-submodule $\Waki^0_{\lambda,k}$, \eqref{def:wakimoto-finite}, which
is isomorphic to the dual Verma module $M^\ast_\lambda$ and satisfies
\eqref{nilp=zero:wakimoto} and \eqref{casimir-on-wakimoto}.
As is mentioned at the end of
\secref{N-point-function}, the restriction of $\Psi(z;q;H)$ to
$\bigotimes_{i=1}^N \Waki^0_{\lambda_i,k}$ satisfies the KZB equation. The
simple structure of $\Waki^0_{\lambda,k}$ makes it possible to write down
the restriction of $\psi(z;t;q;H)$ explicitly.

Each vector $v_i$ in $\Waki^0_{\lambda_i,k}$ corresponds to a polynomial
$P_i(x) \in \Complex[x]$ ($x=(x^\alpha)_{\alpha\in\Delta_+}$) and to an
operator in $\calP_\lambda$ \eqref{def:P-lambda} as
\begin{equation}
    v_i = P_i(\gamma[0]) \ghvac\tensor|\lambda_i\bosket, \qquad
    \Phi_{\lambda_i}(v_i;z) = P_i(\gamma(z))V(\lambda_i;z).
\label{vi=Pi=PiV}
\end{equation}
See \eqref{state-operator-hom:finite}. Let us compute $\psi(z;t;q;H)$ for
this $v_i$. Inserting the expression \eqref{vi=Pi=PiV} into the definition 
\eqref{def:psi}, we obtain
\begin{multline}
    \psi(z;t;q;H)(v_1\tensor\cdots\tensor v_N)
\\
    = q^{-\cvir/24}
    \Tr_{\bosFock_\mu}(
     V(\lambda_1    ;z_1) \cdots V(\lambda_N    ;z_N) 
     V(-\alpha_{i_1};t_1) \cdots V(-\alpha_{i_M};t_M)
    q^{T^\phi[0]} e^{\phi[H;0]})\times
\\
    \times
    \Tr_{\ghFock}(
     P_1(\gamma(z_1)) \cdots P_N(\gamma(z_N))
     \Scr_{\alpha_{i_1}}(t_1) \cdots \Scr_{\alpha_{i_N}}(t_N)
    q^{T^\gh[0]} e^{H^\gh[0]})
    dt_1 \wedge \cdots \wedge dt_M,
\label{psi:boson*ghost}
\end{multline}
where $T^\phi[0]$, $T^\gh[0]$ and $H^\gh[0]$ are zero mode part of
$T^\phi(z)$ \eqref{def:Tphi(z)}, $T^\gh(z)$ \eqref{def:Tgh(z)} and
$H^\gh(z)$ \eqref{explicit:H(z)}, respectively. Since the right hand side
of \eqref{psi:boson*ghost} splits into the bosonic sector and the ghost
sector, we can calculate each part separately.

The computation of the bosonic sector correlation function reduces to
the following lemma. Denote the one-loop correlation function of any
element $A$ of $\widetilde\Bos(\g)$ (\secref{ghost-boson}) by
\begin{equation}
    \langle A \rangle_{\mu,q,H}^\bos 
    :=
    \Tr_{\bosFock_\mu}(A q^{T^\phi[0]}e^{\phi[H;0]}),
\label{def:boson-1-loop}
\end{equation}
when $A|\mu\bosket \in \bosFock_\mu$.
\begin{lem}
\label{lem:boson-1-loop}
Let $\mu_i$ ($i=1\dots,N$) be weights in $\h$ satisfying $\sum_{i=1}^N
\mu_i = 0$. Then the one-loop correlation function of bosonic vertex
operators \eqref{def:boson-vo} is
\begin{multline}
    \langle V(\mu_1;z_1) \cdots V(\mu_N;z_N) 
    \rangle_{\mu,q,H}^{\bos}
    =
    \ell_{\mu_1,\dots,\mu_N,\mu}(z_1,\dots,z_N;q;H)
    :=
    \\ :=
    (q;q)_\infty^{-l} 
    q^{\Delta_\mu}
    e^{(H|\mu)}
    \left( 
      \prod_{i=1}^N
      (\sqrt{-1}\eta(q)^3)^{(\mu_i|\mu_i)/2\kappa} \,
      z_i^{(\mu_i|2\mu - \mu_i)/2\kappa}
    \right)
    \left( 
      \prod_{1\leq i < j \leq N}
      \theta_{11}(z_i/z_j;q)^{(\mu_i|\mu_j)/\kappa},
    \right)
\label{boson-1-loop:formula}
\end{multline}
where $\Delta_\mu=(\mu|\mu+2\rho)/2\kappa$ is the conformal weight and
$\eta(q) = q^{1/24}(q;q)_\infty$ is the Dedekind eta function.
\end{lem}

The proof is in \appref{coherent-states}. This is shown by the
standard method of coherent states (cf.\ for example,
\cite{g-s-w:87}).

The ghost sector can be computed in a similar way as Proposition 3.2 of
\cite{a-t-y:91} and Theorem I of \cite{awa:91}. Let us define the ghost
sector one-loop correlation function by
\begin{equation}
    \langle A \rangle_{q,H}^\gh
    :=
    \Tr_{\ghFock}(A q^{T^\gh[0]}e^{H^\gh[0]}),
\label{def:ghost-1-loop}
\end{equation}
for $A \in \widehat\Gh(\g)$.

The important lemma is the following {\em screening current Ward identity}.

\begin{lem}
\label{scr-Ward}
For any $P_a(x) \in \Complex[x]$ ($a=1,\ldots,n$), a root $\alpha$ and
a sequence of positive roots $\{\alpha(j)\}_{j=1}^m$, we have
\begin{multline}
    \langle
     P_1(\gamma(z_1)) \cdots P_n(\gamma(z_n))
     \Scr_{\alpha}(t)
     \Scr_{\alpha(1)}(t_1) \cdots \Scr_{\alpha(m)}(t_m)
    \rangle_{q,H}^{\gh}
    = \\ =
    \sum_{a=1}^n
    (-w_{\alpha(H)}(t,z_a))
    \langle
    P_1(\gamma(z_1)) \cdots (\Scr_{\alpha}P_a)(\gamma(z_a)) \cdots 
    P_n(\gamma(z_n))
    \Scr_{\alpha(1)}(t_1) \cdots \Scr_{\alpha(m)}(t_m)
    \rangle_{q,H}^\gh
    + \\ +
    \sum_{j=1}^m 
    (-w_{\alpha(H)}(t,t_j))
    f^{\alpha+\alpha(j)}_{\alpha,\alpha(j)}
    \langle
    P_1(\gamma(z_1)) \cdots P_n(\gamma(z_n))
    \Scr_{\alpha(1)}(t_1) \cdots 
    \Scr_{\alpha+\alpha(j)}(t_j) \cdots \Scr_{\alpha(m)}(t_m)
    \rangle_{q,H}^\gh.
\label{scr-Ward:formula}
\end{multline}
Here we used the notations in \eqref{ope:Scr-Scr} and
\eqref{ope:Scr-P}. 
\end{lem}

\begin{proof}
We may assume that $|q|<|t_m|<\cdots<|t_1|<|t|<|z_n|<\cdots<|z_1|<1$. 
The left hand side of \eqref{scr-Ward:formula} is rewritten as follows
because of \eqref{w:pole}:
\begin{equation}
\begin{split}
    &\langle
     P_1(\gamma(z_1)) \cdots P_n(\gamma(z_n))
     \Scr_{\alpha}(t)
     \Scr_{\alpha(1)}(t_1) \cdots \Scr_{\alpha(m)}(t_m)
    \rangle_{q,H}^{\gh}
\\
    =&
    \frac{1}{2\pi i}
    \oint_{\zeta=t} d\zeta\,
    w_{\alpha(H)}(t,\zeta)
    \langle
     P_1(\gamma(z_1)) \cdots P_n(\gamma(z_n))
     \Scr_{\alpha}(\zeta)
     \Scr_{\alpha(1)}(t_1) \cdots \Scr_{\alpha(m)}(t_m)
    \rangle_{q,H}^{\gh}
\\
    =&
    \frac{1}{2\pi i}
    \left(\oint_{|\zeta|=1}
    - \sum_{a=1}^n \oint_{\zeta=z_a}
    - \sum_{j=1}^m \oint_{\zeta=t_j}
    - \oint_{|\zeta|=|q|}
    \right)
    d\zeta \,
    w_{\alpha(H)}(t,\zeta) \times
\\
    &\times
    \langle
     P_1(\gamma(z_1)) \cdots P_n(\gamma(z_n))
     \Scr_{\alpha}(\zeta)
     \Scr_{\alpha(1)}(t_1) \cdots \Scr_{\alpha(m)}(t_m)
    \rangle_{q,H}^{\gh}.
\end{split}
\label{scr-Ward:tmp1}
\end{equation}

The first integral in \eqref{scr-Ward:tmp1} is equal to
\begin{equation}
\begin{split}
    &
    \frac{1}{2\pi i}
    \oint_{|\zeta|=1}
    d\zeta \,
    w_{\alpha(H)}(t,\zeta)
    \Tr_{\ghFock}(
     \Scr_{\alpha}(\zeta)
     P_1(\gamma(z_1)) \cdots P_n(\gamma(z_n))
     \Scr_{\alpha(1)}(t_1) \cdots \Scr_{\alpha(m)}(t_m)
    q^{T^\gh[0]} e^{H^\gh[0]})
\\
    =&
    \frac{1}{2\pi i}
    \oint_{|\zeta|=1}
    d\zeta \,
    w_{\alpha(H)}(t,\zeta)
    \Tr_{\ghFock}(
     P_1(\gamma(z_1)) \cdots P_n(\gamma(z_n))
     \Scr_{\alpha(1)}(t_1) \cdots \Scr_{\alpha(m)}(t_m)
    q^{T^\gh[0]} e^{H^\gh[0]}
    \Scr_{\alpha}(\zeta)
    )
\\
    =&
    \frac{1}{2\pi i}
    \oint_{|\zeta|=1}
    d\zeta \,
    q e^{\alpha(H)}
    w_{\alpha(H)}(t,\zeta) \times
    \\ &\times
    \Tr_{\ghFock}(
     P_1(\gamma(z_1)) \cdots P_n(\gamma(z_n))
     \Scr_{\alpha(1)}(t_1) \cdots \Scr_{\alpha(m)}(t_m)
    \Scr_{\alpha}(q \zeta)
    q^{T^\gh[0]} e^{H^\gh[0]}
    )
\\
    =&
    \frac{1}{2\pi i}
    \oint_{|\zeta|=|q|}
    d\zeta \,
    e^{\alpha(H)}
    w_{\alpha(H)}(t,q^{-1}\zeta)
    \langle
     P_1(\gamma(z_1)) \cdots P_n(\gamma(z_n))
     \Scr_{\alpha(1)}(t_1) \cdots \Scr_{\alpha(m)}(t_m)
    \Scr_{\alpha}(\zeta)
    \rangle_{q,H}^{\gh},
\end{split}
\label{scr-Ward:tmp2}
\end{equation}
where we used the following facts derived from \eqref{ope:H-Scr} and
\eqref{ope:em-Scr}:
\begin{equation}
    e^{H[0]} \Scr_\alpha(\zeta) 
    = e^{\alpha(H)} \Scr_\alpha(\zeta) e^{H[0]}, \qquad
    q^{T[0]} \Scr_\alpha(\zeta) 
    = q \Scr_\alpha(q\zeta) q^{T[0]}.
\label{exp(ad(T[0])/exp(ad(H[0]))Scr(z)}
\end{equation}
Therefore the property \eqref{w:period} of the function
$w_{\alpha(H)}(t,\zeta)$ and \eqref{scr-Ward:tmp2} imply that the
first integral and the last integral in \eqref{scr-Ward:tmp1} cancel.

Using the operator product expansions \eqref{ope:Scr-Scr} and
\eqref{ope:Scr-P}, the second and the third integrals in
\eqref{scr-Ward:tmp1} are rewritten as
\begin{align*}
    &\frac{1}{2\pi i}
    \oint_{\zeta=z_a}
    d\zeta \,
    w_{\alpha(H)}(t,\zeta)
    \langle
     P_1(\gamma(z_1)) \cdots P_n(\gamma(z_n))
     \Scr_{\alpha}(\zeta)
     \Scr_{\alpha(1)}(t_1) \cdots \Scr_{\alpha(m)}(t_m)
    \rangle_{q,H}^{\gh}
\\
    =&
    w_{\alpha(H)}(t,z_a)
    \langle
    P_1(\gamma(z_1)) \cdots (\Scr_{\alpha}P_a)(\gamma(z_a)) \cdots 
    P_n(\gamma(z_n))
    \Scr_{\alpha(1)}(t_1) \cdots \Scr_{\alpha(m)}(t_m)
    \rangle_{q,H}^\gh,
\\
    &\frac{1}{2\pi i}
    \oint_{\zeta=t_j}
    d\zeta \,
    w_{\alpha(H)}(t,\zeta)
    \langle
     P_1(\gamma(z_1)) \cdots P_n(\gamma(z_n))
     \Scr_{\alpha}(\zeta)
     \Scr_{\alpha(1)}(t_1) \cdots \Scr_{\alpha(m)}(t_m)
    \rangle_{q,H}^{\gh}
\\
    =&
    w_{\alpha(H)}(t,t_j) f^{\alpha+\alpha(j)}_{\alpha,\alpha(j)}
    \langle
    P_1(\gamma(z_1)) \cdots P_n(\gamma(z_n))
    \Scr_{\alpha(1)}(t_1) \cdots 
    \Scr_{\alpha+\alpha(j)}(t_j) \cdots 
    \Scr_{\alpha(m)}(t_m)
    \rangle_{q,H}^\gh,
\end{align*}
which proves the lemma.
\end{proof}

\begin{lem}
\label{constant-term}
For any polynomial $P(x) \in \Complex[x]$ with constant term $c_P$,
$$
    \langle P(\gamma(z)) \rangle^\gh_{q,H} = c_P \ch_{\ghFock}(q,H),
$$
where $\ch_{\ghFock}(q,H) = \langle 1 \rangle^\gh_{q,H}$ is the character
of the ghost Fock space. An explicit expression of the character is 
\begin{equation}
    \ch_{\ghFock}(q,H) =
    \prod_{\alpha\in\Delta_+} 
    (e^{-\alpha(H)};q)_\infty^{-1} (q e^{\alpha(H)};q)_\infty^{-1} 
\label{ghost-character}
\end{equation}
\end{lem}

\begin{proof}
It is sufficient to prove that
\begin{equation}
    \Tr_{\ghFock}(\prod_{i=1}^n \gamma_{\beta_i}[n_i]
    q^{T^\gh[0]}e^{H^\gh[0]})
    = 0,
\end{equation}
for any $n\in\Integer_{>0}$, $\beta_i\in\Delta_+$, $n_i\in\Integer$
($i=1,\dots,n$). The Fock space $\ghFock$ has a basis consisting of
vectors of the form
\begin{equation}
    \prod \gamma_{\alpha(i)}[-m_i] 
    \prod \beta_{\alpha'(j)}[-m'_j] \ghvac,
\label{basis:ghFock}
\end{equation}
where $\alpha(i),\alpha'(j)\in\Delta_+$, $m_i\in\Integer_{\geq0}$,
$m'_j\in\Integer_{>0}$. The action of $T^\gh[0]$ and $H^\gh[0]$ is
diagonal with respect to this basis. Hence what we must show is that the
action of $\prod_{i=1}^n \gamma_{\beta_i}[n_i]$ does not have diagonal
components with respect to this basis. This can be shown by elementary
method which uses only the commutation relation \eqref{ghost-relation}.

The character $\ch_{\ghFock}(q;H)$ is calculated by factorizing the total
Fock space into the Fock space $\Fock_{\alpha,m}$ generated by
$\beta_\alpha[m]$ and $\gamma^\alpha[-m]$:
\begin{equation}
    \ghFock = 
    \bigotimes_{\alpha\in\Delta_+,m\in\Integer} \Fock_{\alpha,m}.
\label{ghFock:factorize}
\end{equation}
When $m\geqq 0$, $\Fock_{\alpha,m} = \Complex[\gamma^{\alpha}[-m]]\ghvac$, 
and when $m<0$, $\Fock_{\alpha,m} =
\Complex[\beta^{\alpha}[m]]\ghvac$. The character of each space is:
\begin{equation}
    \Tr_{\Fock_{\alpha,m}}(q^{T^\gh[0]} e^{H^\gh[0]})
    =
    \begin{cases}
    (1-q^m    e^{-\alpha(H)})^{-1},& m\geqq 0,
    \\
    (1-q^{-m} e^{ \alpha(H)})^{-1},& m < 0,
    \end{cases}
\label{ghost-character:alpha,m}
\end{equation}
which follows from the commutation relations,
\begin{alignat*}{2}
    [T^\gh[0],\beta_\alpha[m]]  &= -m \beta_\alpha[m], &\qquad
    [T^\gh[0],\gamma^\alpha[m]] &= -m \gamma^\alpha[m],
\\
    [H^\gh[0],\beta_\alpha[m]]  &=  \alpha(H) \beta_\alpha[m], &\qquad
    [H^\gh[0],\gamma^\alpha[m]] &= -\alpha(H) \gamma^\alpha[m],
\end{alignat*}
and $T^\gh[0]\ghvac = H^\gh[0]\ghvac = 0$.
Multiplying \eqref{ghost-character:alpha,m} over all $\alpha$ and $m$, we
obtain \eqref{ghost-character}.
\end{proof}

\begin{cor}
\label{shapovalov}
For $P(x)\in\Complex[x]$ and a sequence of simple roots
$\{\alpha_{i(j)}\}_{j=1}^n$, we have
$$
    \langle 
      (\Scr_{\alpha_{i(1)}} \cdots \Scr_{\alpha_{i(n)}}P)(\gamma(z))
    \rangle^\gh_{q,H} 
    = \ch_{\ghFock}(q,H)
      (-1)^n \jmath(E_{i(n)} \cdots E_{i(1)} P),
$$
where $\jmath: M^\ast_\lambda \simeqq \Waki^{0}_{\lambda,k} \to \Complex$ is
the pairing with the highest weight vector of the Verma module of $\g$ with
the highest weight $\lambda$, $M_\lambda$. (See
\eqref{def:wakimoto-finite}.) Explicitly written,
\begin{equation}
    \jmath(E_{i(n)} \cdots E_{i(1)} P) =
    R(E_{i(n)})\cdots R(E_{i(1)})P(x)|_{x=0},
\label{jmath}
\end{equation}
where $R(E_i)$ is the differential operator corresponding to the
Chevalley generator $E_i$ given by \eqref{explicit-diff-op}. In
particular, $\jmath(E_{i(n)} \cdots E_{i(1)} P)$ does not depend
on $\lambda$.
\end{cor}

\begin{proof}
According to Lemma 3.3 of \cite{a-t-y:91}, the constant term of
$(\Scr_{\alpha_{i(1)}} \cdots \Scr_{\alpha_{i(n)}}P)(z)$ is given by
$(-1)^n \jmath( E_{i(n)} \cdots E_{i(1)} P) = (-1)^n R(E_{i(n)})\cdots
R(E_{i(1)})P(x)|_{x=0}$.
\end{proof}

\begin{lem}
\label{factorize}
For any $\alpha(i) \in \Delta_+$ ($i=1,\ldots,m$) and
$P_a(x)\in\Complex[x]$ ($a=1,\dots,n$), we have
\begin{multline}
    \frac{
    \langle
      P_1(\gamma(z_1)) \cdots P_n(\gamma(z_n))
      \Scr_{\alpha(1)}(t_1) \cdots \Scr_{\alpha(m)}(t_m)
    \rangle^{\gh}_{q,H}
    }
    {\ch_{\ghFock}(q,H)}
    =\\=
    \sum_{I_1 \sqcup \cdots \sqcup I_n = \{1,\ldots,m\}}
    \prod_{a=1}^n
    \frac{
    \langle
      P_a(\gamma(z_a)) \prod_{i\in I_a} \Scr_{\alpha(i)}(t_i)
    \rangle^\gh_{q,H}
    }
    {\ch_{\ghFock}(q,H)}
  \label{factorize:formula}
\end{multline}
\end{lem}

\begin{proof}
This is a purely combinatorial lemma. We can apply the inductive proof
of (5.3) in \cite{awa:91}, replacing the screening current Ward
identity for genus $0$ with that for genus $1$,
\eqref{scr-Ward:formula}. The first step of the induction (the case
$m=0$) is assured by \lemref{constant-term}.
\end{proof}

\begin{lem}
\label{SP}
For any $P(x)\in\Complex[x]$ and roots $\alpha(i)$
($i=1,\dots,m$), we have
\begin{multline}
    \langle 
      P(\gamma(z)) \Scr_{\alpha(1)}(t_1) \cdots \Scr_{\alpha(m)}(t_m)
    \rangle^\gh_{q,H} = \\ =
    \sum_{\sigma\in\germanS_m}
    (-w_{\alpha(\sigma(1))}(t_{\sigma(1)}, t_{\sigma(2)}))
    (-w_{\alpha(\sigma(1))+\alpha(\sigma(2))}
    (t_{\sigma(2)}, t_{\sigma(3)})) \cdots 
    (-w_{\alpha(\sigma(1))+ \cdots + \alpha(\sigma(m))}
    (t_{\sigma(m)}, z)) \times \\ \times 
    \langle
     (\Scr_{\alpha(\sigma(1))} \cdots
      \Scr_{\alpha(\sigma(m))}P)(\gamma(z)) 
    \rangle^\gh_{q,H},
\label{SP:formula}
\end{multline}
where we write $w_{\alpha(H)}$ as $w_{\alpha}$ for short.
\end{lem}

\begin{proof}
We prove this statement by induction on $m$, as in the proof of (5.4)
in \cite{awa:91}. When $m=0$, the statement is trivial and when $m=1$, 
it is nothing but the screening current Ward identity
\eqref{scr-Ward:formula}.

Assume that \eqref{SP:formula} holds for all $m \leqq n$. Let us
regard the left and right hand side of \eqref{SP:formula} for $m=n+1$
as functions of $t_0$:
\begin{align}
    F_1(t_0) :=& 
    \langle 
      P(\gamma(z)) 
      \Scr_{\alpha(0)}(t_0) 
      \Scr_{\alpha(1)}(t_1) \cdots \Scr_{\alpha(n)}(t_n)
    \rangle^\gh_{q,H},
\label{def:F1}
\\
    F_2(t_0) :=&
    \sum_{\sigma\in\germanS_{n+1}} 
    (-w_{\alpha(\sigma(0))}(t_{\sigma(0)}, t_{\sigma(1)}))
    (-w_{\alpha(\sigma(0))+\alpha(\sigma(1))}
    (t_{\sigma(1)}, t_{\sigma(2)})) \cdots 
    (-w_{\alpha(\sigma(0))+ \cdots + \alpha(\sigma(n))}
    (t_{\sigma(n)}, z)) \times \nonumber \\ &\times 
    \langle
     (\Scr_{\alpha(\sigma(1))} \cdots
      \Scr_{\alpha(\sigma(n))}P)(\gamma(z)) 
    \rangle^\gh_{q,H},
\label{def:F2}
\end{align}
where $\sigma$ is a permutation,
$\sigma:\{0,1,\dots,n\}\to\{0,1,\dots,n\}$. We now show that
$F_1(t_0)=F_2(t_0)$.

First, note that both functions are meromorphic on $\Complex^\times$
and poles exist at $t_0=t_i$ ($i=1,\dots,n$) and at $t_0=z$.

(i) Both functions have the same quasi-periodicity,
\begin{equation}
    f(qt_0) = e^{-\alpha(0)(H)}q^{-1} f(t_0).
\label{f:period}
\end{equation}
In fact, \eqref{f:period} is proved for $f=F_1$ similarly to
\eqref{scr-Ward:tmp2}. It follows from the property of the function
$w$ \eqref{w:period} that $F_2(t_0)$ also satisfies the same
periodicity property \eqref{f:period}.

(ii) The principal parts of the pole at $t_0=z$ are equal to
\begin{equation}
    \frac{1}{t_0-z}
    \langle
     (\Scr_{\alpha(0)}P)(\gamma(z)) 
     \Scr_{\alpha(1)}(t_1) \cdots \Scr_{\alpha(n)}(t_n)
    \rangle^\gh_{q,H}.
\label{pole:t0=z}
\end{equation}
For $F_1(t_0)$, this is a direct consequence of the Ward identity
\eqref{scr-Ward:formula}. The pole of $F_2(t_0)$ at $t_0=z$ comes from
terms in \eqref{def:F2} such that $\sigma(n)=0$. Using \eqref{w:pole}
and the induction hypothesis, we can show that its principal part is
of the form \eqref{pole:t0=z}.

(iii) The principal parts of the pole at $t_0=t_i$ are equal to
\begin{equation}
    \frac{f^{\alpha(0)+\alpha(i)}_{\alpha(0),\alpha(i)}}{t_0-t_i}
    \langle
     P(\gamma(z))
     \Scr_{\alpha(1)}(t_1) \cdots \Scr_{\alpha(0)+\alpha(i)}(t_i) \cdots
     \Scr_{\alpha(n)}(t_n)
    \rangle^\gh_{q,H}.
\label{pole:t0=ti}
\end{equation}
The Ward identity \eqref{scr-Ward:formula} implies \eqref{pole:t0=ti} for
$F_1(t_0)$. The pole of $F_2(t_0)$ at $t_0=t_i$ comes from terms in
\eqref{def:F2} such that $(\sigma(0),\sigma(i))=(j-1,j)$ or
$(\sigma(0),\sigma(i))=(j,j-1)$ ($j=1,\dots,n$). The principal part
becomes
\begin{multline}
    \frac{1}{t_0-t_i}
    \sum_{j=1}^n \sum_{\sigma}
    (-w_{\alpha(\sigma(0))}(t_{\sigma(0)},t_{\sigma(1)}))
    \cdots
\\
    \cdots
    (-w_{\alpha(\sigma(0))+\cdots+\alpha(\sigma(j-2))}
    (t_{\sigma(j-2)},t_i))
    (-w_{\alpha(\sigma(0))+\cdots+\alpha(\sigma(j-2))+\alpha(0)+\alpha(i)}
    (t_i, t_{\sigma(j+1)}))
    \cdots
\\
    \langle
     (\Scr_{\alpha(\sigma(0))}\cdots [\Scr_{\alpha(0)},\Scr_{\alpha(i)}]
     \cdots
     \Scr_{\alpha(\sigma(n))}P)(z)
    \rangle^\gh_{q,H}.
\label{pole:t0=ti:F2}
\end{multline}
where $\sigma$ runs through the set of permutations
$\sigma:\{0,\dots,j-2,j+1,\dots,n\}\to\{1,\dots,i-1,i+1,\dots,n\}$. Since
$[\Scr_{\alpha(0)},\Scr_{\alpha(i)}]=
f^{\alpha(0)+\alpha(i)}_{\alpha(0),\alpha(i)}\Scr_{\alpha(0)+\alpha(i)}$,
it follows from the induction hypothesis that \eqref{pole:t0=ti:F2} is
equal to \eqref{pole:t0=ti}.

Comparing $F_1(t_0)$ and $F_2(t_0)$ by (i), (ii) and (iii), we conclude
that $F_1(t_0)=F_2(t_0)$.
\end{proof}

Putting together \eqref{psi:boson*ghost}, \lemref{lem:boson-1-loop},
\lemref{factorize}, \lemref{SP}, \corref{shapovalov}, and
\eqref{ghost-character}, we finally obtain the integral representation of
a solution of the KZB equations.
\begin{thm}
\label{integral-rep}
The following integral gives a solution of the KZB equations,
\eqref{kzb:weight0}, \eqref{kzb} with $c_2^{(j)} =
(\lambda_j|\lambda_j+2\rho)$ and \eqref{kzb-heat}:
\begin{equation}
    \tilde\Psi^0(z;q;H) =
    \int_{{\mathcal C}(z,q,H)}
    \ell^0_{-\alpha_{i_1},\dots,-\alpha_{i_M},\lambda_1,\dots,\lambda_N,\mu}
    (t;z;q;H)
    \psi^\gh(t;z;q;H;P_1,\dots,P_N),
\label{integral-rep:formula}
\end{equation}
where ${\mathcal C}(z,q,H)$ is a family of $M$-cycles with coefficients in
$\Line^\ast$,
\begin{multline}
    \ell^0_{\mu_1,\dots,\mu_N,\mu}(z_1,\dots,z_N;q;H)
    :=
    \\ :=
    q^{(\mu+\rho|\mu+\rho)/2\kappa}
    e^{\mu(H)}
    \left( 
      \prod_{i=1}^N
      (\sqrt{-1}\eta(q)^3)^{(\mu_i|\mu_i)/2\kappa} \,
      z_i^{(\mu_i|2\mu - \mu_i)/2\kappa}
    \right)
    \left( 
      \prod_{1\leq i < j \leq N}
      \theta_{11}(z_i/z_j;q)^{(\mu_i|\mu_j)/\kappa},
    \right)
\label{def:ell0}
\end{multline}
and the $M$-form $\psi^\gh$ is defined as follows: $P_a(x)$ are
polynomials in $x$,
\begin{equation}
    \psi^\gh(t;z;q;H;P_1,\dots,P_N)
    =
    e^{\rho(H)}
    \sum_{I_1 \sqcup \cdots \sqcup I_N = \{1,\ldots,M\}}
    \prod_{a=1}^N
    \frac{
    \langle
      P_a(\gamma(z_a)) \prod_{j\in I_a} \Scr_{\alpha_{i_j}}(t_j)
    \rangle^\gh_{q,H}
    }{\ch_{\ghFock}(q,H)}
    dt_1 \wedge \cdots \wedge dt_M,
\label{def:psi-gh}
\end{equation}
and the last factor in \eqref{def:psi-gh} for $a$ ($1\leqq a \leqq N$) is
\begin{multline}
    \frac{
    \langle 
      P(\gamma(z)) 
      \Scr_{\alpha_{i(1)}}(t_1) \cdots \Scr_{\alpha_{i(m)}}(t_m)
    \rangle^\gh_{q,H}}
    {\ch_{\ghFock}(q,H)}
    = \\ =
    \sum_{\sigma\in\germanS_m} 
    w_{\alpha_{i(\sigma(1))}}(t_{\sigma(1)}, t_{\sigma(2)})
    w_{\alpha_{i(\sigma(1))}+\alpha_{i(\sigma(2))}}
    (t_{\sigma(2)}, t_{\sigma(3)}) \cdots 
    w_{\alpha_{i(\sigma(1))}+ \cdots + \alpha_{i(\sigma(m))}}
    (t_{\sigma(m)}, z) \times \\ \times
    \jmath(E_{i(\sigma(m))} \cdots E_{i(\sigma(1))}P),
\end{multline}
if $\{i_j\,|\,j\in I_a\}=\{i(1),\dots,i(m)\}$.
\end{thm}
This result is an elliptic analogue of \cite{sch-var:90},
\cite{sch-var:91}, \cite{a-t-y:91}, \cite{awa:91} and a generalization of
a result for $sl(2)$ in \cite{ber-fel:90}. Felder and Varchenko
have obtained a similar formula in \cite{fel-var:95} from a different
standpoint.

\section{Concluding remarks}
\label{conclusion}

We found an integral representation of $N$-point functions of the WZW
model on elliptic curves and gives an explicit expression for a solution
of the KZB equations, using the Wakimoto realization. Let us list some of
related problems.

\begin{enumerate}
\item
Higher genus: Is there a similar integral representations of correlation
functions of the Wess-Zumino-Witten models on higher genus Riemann
surfaces? There are several works to this direction \cite{gmmos:90},
\cite{kon:92}. Their formulations are, however, different from ours. 

\item
Twisted Wess-Zumino-Witten models: In \cite{kur-tak:97} we formulated
``another'' Wess-Zumino-Witten model on elliptic curves which we named a
``twisted WZW model''. Is it possible to give an integral representation
of solutions of the KZ type equations for the correlation functions?

\item
Critical level: Feigin, Frenkel and Reshetikhin \cite{f-f-r:94} found that
the Bethe vector of a certain spin chain model is obtained from the
Wakimoto realization of the Wess-Zumino-Witten model on the Riemann sphere
at the critical level. How about the genus one case?
\end{enumerate}

We shall study the last question in the forthcoming paper. In fact, Felder
and Varchenko \cite{fel-var:95} have found a relation of their integral
representation of solutions of the KZB equations with a solution of
a quantum $N$-body system. We shall take the conformal field theoretical
approach to this problem. 

\subsection*{Acknowledgments}
TT is supported by grant-in-aid of the Ministry of Education and Sciences
of Japan, No.~09740009. The authors express their gratitude to
Edward~Frenkel, Takeshi~Ikeda, Akishi~Kato, Hitoshi~Konno, Atsushi~Matsuo,
Alexei~Morozov, Takeshi~Suzuki and Yasuhiko~Yamada for comments and
discussions.

\appendix
\section{Theta functions}
\label{theta}

We denote the theta function with characteristic $(1/2, 1/2)$ (cf.\ Chapter 
I of \cite{mum}) additively as
\begin{equation}
    \theta_{11}(x;\tau) = 
    \sum_{n\in\Integer} e^{(n+1/2)^2 \pi i \tau + 2\pi i (n+1/2)(x+1/2)}
\label{def:theta:add}
\end{equation}
and multiplicatively as
\begin{equation}
    \theta_{11}(z;q) = \theta_{11}(x;\tau),
\label{def:theta:mult}
\end{equation}
where $z=\exp(2\pi i x)$ and $q=\exp(2\pi i \tau)$. The infinite product
expansion
\begin{equation}
    \theta_{11}(z;q) 
    =
    i (q;q)_\infty q^{1/8} z^{1/2} (z^{-1};q)_\infty (qz;q)_\infty
\label{theta:inf-prod}
\end{equation}
is also useful, where $(x;q)_\infty = \prod_{n=0}^\infty (1-x q^n)$.

We use a function $w_c(w,z)$ on $\Complex^\times \times \Complex^\times$
with parameter $c \in \Complex^\times$ characterized by the following
properties:
\begin{enumerate}
\item
$w_c(w,z)$ is a meromorphic function of $z$ and $w$.

\item
$w_c(w,z)$ has a following (quasi-)periodicity
\begin{equation}
    w_c(w,qz) = e^c w_c(w,z), \qquad
    w_c(qw,z) = q^{-1} e^{-c} w_c(w,z).
\label{w:period}
\end{equation}

\item
$w_c(w,z)$ has only one simple pole on the elliptic curve
$\Complex^\times/q^\Integer$ at $z = w$ as a function of $z$. Its Laurent
expansion around $z=w$ is:
\begin{equation}
    w_c(w,z) = \frac{1}{z-w} + \text{regular}.
\label{w:pole}
\end{equation}
\end{enumerate}

An explicit form of $w_c(w,z)$ is as follows:
\begin{equation}
    w_c(w,z) = 
    \frac{\theta'_{11}(1;q)}{\theta_{11}(e^{-c};q)}
    \frac{\theta_{11}(e^{-c} z/w;q)}{w \theta_{11}(z/w;q)},
\end{equation}
where $\theta'_{11}(z;q)=d/dz \theta_{11}(z;q)$.

To write down the Knizhnik-Zamolodchikov-Bernard equations, we need
following functions:
\begin{align}
    \sigma_c(z) &:= 
    \frac{\theta'_{11}(1;q)}{\theta_{11}(e^{-c};q)}
    \frac{\theta_{11}(e^{-c} z;q)}{\theta_{11}(z;q)},
\label{def:sigma}
\\
    \zeta(z) &:= z \frac{\theta'_{11}(z;q)}{\theta_{11}(z;q)}.
\label{def:zeta}
\end{align}

\section{Method of coherent states and one-loop correlation functions}
\label{coherent-states}

In this appendix we review the method of coherent states, following
Chapter 7.A and 8.1 of \cite{g-s-w:87} and compute one-loop correlation
functions of vertex operators of the free bosons, which proves
\lemref{lem:boson-1-loop}. 

Let $a$ and $a^\dagger$ be generators of a Heisenberg algebra $\Heis$:
\begin{equation}
    [a, a^\dagger] = 1,
\label{com-rel:a}
\end{equation}
and $|0\rangle$ and $\langle 0|$ be the generating vector of the Fock
space representation of $\Heis$ and that of its (restricted) dual,
respectively:
\begin{align}
    \Fock &:= \Heis |0\rangle, \qquad  a |0\rangle = 0,
\label{def:fock}
\\
    \Fock^\ast &:= \langle 0|\Heis, \qquad  \langle 0| a^\dagger = 0.
\label{def:fock*}
\end{align}
There are natural bases $\{|n\rangle\}_{n\in\Natural}$ of $\Fock$ and
$\{\langle n|\}_{n\in\Natural}$ of $\Fock^\ast$, consisting of
eigenvectors of the number counting operator $N_a=a^\dagger a$:
\begin{equation}
\begin{gathered}
    |n\rangle = \frac{(a^\dagger)^n}{\sqrt{n!}}|0\rangle, \qquad
    \langle n| = \langle 0| \frac{a^n}{\sqrt{n!}},
\\
    N_a |n \rangle = n |n \rangle, \qquad
    \langle n| N_a = \langle n| n, \qquad
    \langle m|n \rangle = \delta_{mn}.
\end{gathered}
\label{number-basis}
\end{equation}
The {\em coherent states} are defined by
\begin{equation}
\begin{aligned}
    |\lambda) &:= \exp(\lambda a^\dagger) |0\rangle
    = \sum_{n=0}^\infty \frac{\lambda^n}{\sqrt{n!}} |n\rangle,
\\
    (\lambda| &:= \langle0|\exp(\bar{\lambda} a)
    = \sum_{n=0}^\infty \langle n| \frac{\bar{\lambda}^n}{\sqrt{n!}},
\end{aligned}
\label{def:coherent-states}
\end{equation}
for $\lambda\in\Complex$. Here $\bar{\lambda}$ is the complex conjugate of
$\lambda$. In particular, $|0) = |0\rangle$ and $(0| = \langle 0|$. They
have the following properties:
\begin{gather}
    a |\lambda) = \lambda |\lambda), \qquad
    (\lambda| a^\dagger = (\lambda| \bar{\lambda},
\label{eigenvector-of-ann/cr-op}
\\
    (\mu|\lambda) = e^{\bar{\mu}\lambda},
\label{coherent-states:inner-prod}
\\
    q^{N_a} |\lambda) = |q\lambda),
\label{coherent-states:shift}
\end{gather}
where $q\in\Complex^\times$. The trace of an operator
$A\in\End_\Complex(\Fock)$ is computed by the following integral:
\begin{equation}
    \Tr_{\Fock}(A) =
    \frac{1}{\pi}\int_{\Complex} d^2 \lambda \,
    e^{-|\lambda|^2} (\lambda|A|\lambda).
\label{trace-by-coherent-states}
\end{equation}

Using these formulae, we calculate the one-loop correlation function
\eqref{def:boson-1-loop} of vertex operators of the free boson fields
$V(\mu_i;z)$ \eqref{def:boson-vo}:
\begin{equation}
    \langle V(\mu_1;z_1) \cdots V(\mu_N;z_N)
    \rangle_{\mu,q,H}^\bos 
    =
    \Tr_{\bosFock_\mu}
    (V(\mu_1;z_1) \cdots V(\mu_N;z_N) q^{T^\phi[0]}e^{\phi[H;0]}), 
\label{boson-vo-1-loop}
\end{equation}
where $\sum_{i=1}^N \mu_i = 0$. 

Let us fix an orthonormal basis $\{H_r\}_{r=1}^l$ of $\h$. The boson Fock
space $\bosFock_\mu$ is factorized as
\begin{equation}
    \bosFock_\mu = 
    \Fock_{0,\mu} \tensor 
    \bigotimes_{r=1}^l \bigotimes_{n=1}^\infty \Fock_{r,n},
\label{boson-fock-factor}
\end{equation}
where $\Fock_{0,\mu}$ is the zero-mode space $\Complex |\mu\bosket$
and $\Fock_{r,n}$ is a non-zero-mode Fock space generated by
$\phi_r[\pm m]:=\phi[H_r;\pm n]$. Note that $\phi_r[\pm m]$ and
$\phi_{r'}[\pm n]$ ($m,n\in\Integer_{>0}$, $r,r'=1,\dots,l$) commute
with each other unless $r=r'$ and $m=n$. Hence, to compute the value
of \eqref{def:boson-1-loop}, we have to compute the trace of
$V(\mu_1;z_1) \cdots V(\mu_N;z_N) q^{T^\phi[0]}e^{\phi[H;0]}$
over $\Fock_{0,\mu}$ and $\Fock_{r,n}$ and multiply all of them. The
vertex operators $V(\mu_i;z_i)$ are factorized into product of the
zero mode part $V_0(\mu_i;z_i)$ (\eqref{def:V0}) and the non-zero
mode part $\tilde V(\mu_i;z_i)$ (\eqref{def:tildeV}), while the
operator $T^\phi(z)$ is decomposed into sum of the zero mode and the
non-zero mode parts as
\begin{align}
    T^\phi(z) &= T^\phi_0(z) + {\tilde T}^\phi(z),
\label{Tphi(z):decomposition}
\\
    T^\phi_0(z) 
    &:= \frac{1}{2\kappa}\sum_{i=1}^l \phi[H_i;0] \phi[H^i;0]
      + \frac{1}{2\kappa}\phi[2\rho;0], 
\label{def:Tphi0}
\\
    {\tilde T}^\phi(z)
    &:= \frac{1}{\kappa}
    \sum_{r=1}^l \sum_{n=1}^\infty \phi_r[-n] \phi_r[n].
\label{def:tildeTphi}
\end{align}

(I) Zero-mode:
It is easy to see that the zero mode part of $V(\mu_1;z_1) \cdots
V(\mu_N;z_N) q^{T^\phi[0]}e^{\phi[H;0]}$ acts on $|\mu\bosket$ as
\begin{multline}
    V_0(\mu_1;z_1) \cdots V_0(\mu_N;z_N)
    q^{T^\phi_0[0]}e^{\phi[H;0]} 
    |\mu\bosket
\\
    =
    \prod_{1\leq i < j \leq N} z_i^{(\mu_i|\mu_j)/\kappa}
    \prod_{i=1}^N z_i^{(\mu_i|\mu)/\kappa}
    q^{(\mu + 2\rho|\mu)/2\kappa}
    e^{\mu(H)}
    |\mu\bosket.
\label{boson-vo-1-loop:0-mode}
\end{multline}

(II) Non-zero-mode:
The algebra $\Heis_{r,n}$ generated by $\phi_r[\pm n]$ is isomorphic to
$\Heis$ through the isomorphism defined by
$$
    a = \frac{\phi_r[n]}{\sqrt{\kappa n}}, \qquad
    a^\dagger = \frac{\phi_r[-n]}{\sqrt{\kappa n}}.
$$
The $\Heis_{r,n}$ part of $V(\mu_1;z_1) \cdots V(\mu_N;z_N)
q^{T^\phi[0]}e^{\phi[H;0]}$ is
$$
    \exp
    \left(\frac{\mu_1^r}{\kappa} \frac{z_1^n}{n} \phi_r[-n]\right)
    \exp
    \left(\frac{\mu_1^r}{\kappa} \frac{z_1^{-n}}{-n} \phi_r[n]\right) 
    \cdots
    \exp
    \left(\frac{\mu_N^r}{\kappa} \frac{z_N^n}{n} \phi_r[-n]\right)
    \exp
    \left(\frac{\mu_N^r}{\kappa} \frac{z_N^{-n}}{-n} \phi_r[n]\right)
    q^{\phi_r[-n]\phi_r[n]/\kappa},
$$
where $\mu_i = \sum_{r=1}^l \mu_i^r H_r$. Its trace over
$\Fock_{r,n}$ is computed by means of the formula
\eqref{trace-by-coherent-states}. The result is
\begin{equation}
    \frac{1}{1-q^n}
    \exp \left(
    \frac{-1}{\kappa n(1-q^n)}
    \sum_{1\leq i < j \leq N} \mu_i^r \mu_j^r 
        \left(\frac{z_j}{z_i}\right)^n
    +
    \frac{-q^n}{\kappa n(1-q^n)}
    \sum_{1\leq i \leq j \leq N} \mu_i^r \mu_j^r 
        \left(\frac{z_i}{z_j}\right)^n
    \right).
\label{boson-vo-1-loop:(r,n)}
\end{equation}
Multiplying \eqref{boson-vo-1-loop:(r,n)} for all $r$ and $n$, we have 
\begin{multline}
    \Tr_{\bigotimes_{r,n}\Fock_{r,n}}
    (\tilde V(\mu_1;z_1) \cdots \tilde V(\mu_N;z_N)
    q^{{\tilde T}^\phi[0]})
    =\\=
    (q;q)_\infty^{-l}
    \prod_{i=1}^N (q;q)_\infty^{(\mu_i|\mu_i)/\kappa}
    \prod_{1\leq i < j \leq N}
    \bigl(
    (z_j/z_i;q)_\infty (q z_i/z_j;q)_\infty
    \bigr)^{(\mu_i|\mu_j)/\kappa}.
\label{boson-vo-1-loop:non-0-mode}
\end{multline}

Putting \eqref{boson-vo-1-loop:0-mode} and
\eqref{boson-vo-1-loop:non-0-mode} together, we obtain the final
result \eqref{boson-1-loop:formula}, using the infinite product
expansion of the theta function \eqref{theta:inf-prod}. This completes
the proof of \lemref{lem:boson-1-loop}.


\end{document}